\documentclass{amsart}

\usepackage{lineno}
\usepackage{graphicx}
\usepackage{mathrsfs}
\usepackage{stmaryrd}
\usepackage{amsfonts}

\usepackage{enumerate,amsmath,amssymb,amsthm}

\usepackage{geometry}
\renewcommand{\theequation}{\arabic{equation}}
\newcommand{\be}{\begin{eqnarray}}
\newcommand{\ee}{\end{eqnarray}}
\newcommand{\ce}{\begin{eqnarray*}}
\newcommand{\de}{\end{eqnarray*}}
\newtheorem{theorem}{Theorem}[section]
\newtheorem{lemma}[theorem]{Lemma}
\newtheorem{remark}[theorem]{Remark}
\newtheorem{definition}[theorem]{Definition}
\newtheorem{proposition}[theorem]{Proposition}
\newtheorem{Examples}[theorem]{Examples}
\newtheorem{corollary}[theorem]{Corollary}

\def\[{{\Big[}}
\def\]{{\Big]}}
\def\<{{\langle}}
\def\>{{\rangle}}
\def\({{\Big(}}
\def\){{\Big)}}

\def\bx{{\mathbf{x}}}

\def\bt{\begin{theorem}}
\def\et{\end{theorem}}
\def\bl{\begin{lemma}}
\def\el{\end{lemma}}
\def\br{\begin{remark}}
\def\er{\end{remark}}
\def\bx{\begin{Examples}}
\def\ex{\end{Examples}}
\def\bd{\begin{definition}}
\def\ed{\end{definition}}
\def\bp{\begin{proposition}}
\def\ep{\end{proposition}}
\def\bc{\begin{corollary}}
\def\ec{\end{corollary}}

\def\geq{\geqslant}
\def\leq{\leqslant}

\theoremstyle{thmit} 

\theoremstyle{thmrm} 

\newtheorem*{oldproof}{Proof}

\title{Slow manifold and parameter estimation for a nonlocal fast-slow stochastic evolutionary system\footnote{The research   was partly supported  by the NSF   grant 1620449   and NSFC grants  11531006 and  11771449.}}

\author{Hina Zulfiqar}
\address{School of Mathematics and Statistics,
Huazhong University of Science and Technology,\\ Wuhan 430074, China \\ \& Center for  Mathematical Sciences, Huazhong University of Science and Technology\\zhinazulfiqar@gmail.com}
\author{Ziying He}
\address{School of Mathematics and Statistics,
Huazhong University of Science and Technology,\\ Wuhan 430074, China \\ \& Center for  Mathematical Sciences, Huazhong University of Science and Technology\\ziyinghe@hust.edu.cn}
\author{Meihua Yang}
\address{School of Mathematics and Statistics,
Huazhong University of Science and Technology,\\ Wuhan 430074, China \\yangmeih@hust.edu.com}
\author{Jinqiao Duan}
\address{Department of Applied Mathematics, Illinois Institute of Technology,\\ Chicago, IL 60616, USA \\ \& Center for  Mathematical Sciences, Huazhong University of Science and Technology\\duanjq@gmail.com}

\begin{document}

\maketitle
\begin{abstract}
 \indent We establish a slow manifold for a fast-slow stochastic evolutionary system with anomalous diffusion, where both fast and slow components are influenced by white noise. Furthermore, we  prove the exponential tracking property for the random slow manifold and this leads to a lower dimensional reduced system based on the slow manifold. Also we consider parameter estimation for this nonlocal fast-slow stochastic dynamical system, where only the slow component is observable. In quantifying parameters in stochastic evolutionary systems, this offers an advantage of dimension reduction. \\
\textbf{keywords:} Nonlocal Laplacian, fast-slow stochastic system, random slow manifold, exponential tracking property, parameter estimation.
\end{abstract}

\section{Introduction}\label{s:1}
\begin{linenomath*}
 Nonlocal operators arise in complex phenomena such as anomalous diffusion in geophysical flows \cite{caffarelli2010drift,meerschaert2012stochastic,metzler2004restaurant}.   The usual Laplacian operator $\Delta$  is the Markov generator of  the Gaussian process called   Brownian motion (or Wiener process). It   may be regarded as the  macroscopic manifestation of this normal diffusion phenomenon. Thus we often see the Laplacian operator $\Delta$ in models for diffusion-related complex systems appear as stochastic or deterministic   evolutionary differential equations.    An exceptional, yet imperative nonlocal operator is nonlocal Laplacian operator emerged in non-Gaussian or anomalous diffusion. The nonlocal Laplacian operator $(-\Delta)^{\frac{\alpha}{2}}$ is the Markov generator of the  non-Gaussian process  called symmetric $\alpha$-stable L$\acute{e}$vy motion, for $\alpha\in(0,2)$,  and it is the  macroscopic manifestation of this anomalous diffusion phenomenon \cite{applebaum2009levy, duan2015introduction}.  \\
\indent For investigating dynamical behaviours of deterministic systems, the theory of invariant manifolds have a satisfactory and popular history. Hadamard \cite{hadamard1901iteration} introduced
it for the first time, then by Lyapunov and Perron \cite{caraballo2004existence,duan2004smooth,chow1988invariant}. For deterministic system it has been modified by numerous authors \cite{duan2003invariant, ruelle1982characteristic,bates1998existence,chicone1997center,chow1991smooth,henry1981lecture}. Duan and Wang \cite{duan2014effective} explained effective dynamics for the stochastic partial differential equations. Schmalfu${\ss}$ and Schneider \cite{schmalfuss2008invariant} have recently explored random inertial manifolds for stochastic differential equations with two times scales by a fixed point technique, which eliminate the fast variables via random graph transformation. They showed that as the scale parameter approaches to zero the inertial manifold approaches to some other manifold, which can be defined in term of slow manifold. Wang and Roberts \cite{wang2013slow} explored the averaging of slow manifold for fast-slow stochastic differential equations, in which only fast mode is influenced by white noise. For multi-time-scale stochastic system, the exploration of slow manifolds was conducted by Fu, Liu and Duan \cite{fu2012slow}. Bai $\emph{et al.}$  \cite{bai2017slow} recently established the existence of slow manifold for nonlocal fast-slow stochastic evolutionary equations in which only fast mode influenced by white noise. Ren $\emph{et al.}$ \cite{ren2015parameter} devised a parameter estimation method based on a random slow manifold for a finite dimensional slow-fast stochastic system. \\ \indent In the present paper we consider invariant manifold for nonlocal fast-slow stochastic evolutionary system, in which fast and slow mode, both are influenced by white noise. Namely, we explore the following nonlocal fast-slow stochastic evolutionary system.
\begin{align} &\dot{u}^{\epsilon}=-\frac{1}{\epsilon}(-\Delta)^{\frac{\alpha}{2}}u^{\epsilon}+\frac{1}{\epsilon}f(u^{\epsilon},v^{\epsilon})+\frac{\sigma_{1}}{\sqrt{\epsilon}}\dot{W}_{t}^{1},\mbox{ in }H_{1}\\&\dot{v}^{\epsilon}=Jv^{\epsilon}+g(u^{\epsilon},v^{\epsilon})+\sigma_{2}\dot{W}_{t}^{2},\mbox{ in }H_{2}\\ &u^{\epsilon}|(-1,1)^{c}=0,\indent \indent v^{\epsilon}|(-1,1)^{c}=0, \end{align}where, for $x\in \mathbb{R}$ and $\alpha\in(0,2),$\\$$(-\Delta)^{\frac{\alpha}{2}}u^{\epsilon}(x,t)=\frac{2^{\alpha}\Gamma(\frac{1+\alpha}{2})}{\sqrt{\pi}|\Gamma(\frac{-\alpha}{2})|}P.V.\int_{\mathbb{R}}\frac{u^{\epsilon}(x,t)-v^{\epsilon}(y,t)}{|x-y|^{1+\alpha}}dy,$$ \\is called fractional Laplacian operator, which is nonlocal operator with the Cauchy principle value $(P.V.)$ considering as the limit of integral over $\mathbb{R}\setminus(x-\epsilon,x+\epsilon)$ for $\epsilon$ tends to zero. The Gamma function is defined to be\\$$\Gamma(y)=\int_{0}^{\infty}t^{y-1}e^{-t}dt,\indent \forall \indent y>0.$$\indent  The fundamental objective of this paper is to establish the existence of random invariant manifold $M^{\epsilon}$, for small enough $\epsilon>0$ with an exponential tracking property in Section 4, and conduct up a parameter estimation by using observations only on the slow system for the above nonlocal fast-slow stochastic evolutionary system in Section 5. Since slow system is lower dimensional as compared to the original stochastic system, so in computational cost this parameter estimator offers a benefit. In addition, it provides advantage of using the observations only on slow components. We verified the results about slow manifold, exponential tracking property and parameter estimation in last section via theoretical example and a example from numerical simulation. Before constructing the invariant manifold, initially we need to prove the existence, and then uniqueness of solution for newly made system (1)-(2). So it is verified in next Section. Many nonlocal problems explained in \cite{bucur2016nonlocal} are just similar to (1)-(3). There are usually two techniques for the establishment of invariant manifolds: Hadamard graph transform method \cite{schmalfuss1998random,duan2003invariant} and Lyapunov-Perron method \cite{chow1988invariant,duan2004smooth, caraballo2004existence} . For the achievement of our objective we used Lyapunov method, which is different from method of random graph transformation. More specifically, it can be shown that this manifold $M^{\epsilon}$ asymptotically approximated to a slow manifold $M^{0}$ for sufficiently small $\epsilon$, under suitable conditions.\\ \indent The state space for the system $(1)-(3)$ is a separable Hilbert space $\mathbb{H}=H_{1}\times H_{2},$ which is a product of Hilbert spaces $H_{1}$ and $H_{2}$. Norm $||\cdot||_{1}$ and $||\cdot||_{2}$ are considered as the norm of $H_{1}$ and $H_{2}$ respectively. The norm of separable Hilbert space $\mathbb{H}$ is defined by$$||\cdot||=||\cdot||_{1}+||\cdot||_{2}.$$\indent In the system $(1)-(3)$, $\epsilon$ is a parameter with the bound $0<\epsilon\ll1$, which denotes the ratio of two times scales such that $||\frac{dv}{dt}||_{H_{2}}\ll||\frac{du}{dt}||_{H_{1}}.$\\ \indent In order to balance the stochastic force and deterministic force, strength of noise in fast equation is chosen to be $\sigma_{1}/\sqrt{\epsilon}$, while in slow equation is $\sigma_{2}$ respectively. The interaction functions $f$ and $g$ are non linearities. Linear operator $J$ generates a $C_{0}$-semigroup, which satisfy the slow evolution hypothesis given in next section. Consider that $u^{\epsilon}\in H_{1}$ is the fast mode, while $v^{\epsilon} \in H_{2}$ is the slow mode.\\ \indent The Weiner process $W_{t}^{1}$ and $W_{t}^{2}$ are time dependent Brownian motions and defined on the probability space $\Omega_{1}$ and $\Omega_{2}$.  While natural filtrations $\mathcal{F}_{1}$ and $\mathcal{F}_{2}$  are generated by $W^{1}$ and $W^{2}$ respectively. The white noises $\dot{W}_{t}^{1}$ and $\dot{W}_{t}^{2}$ are generalized time derivatives of $W^{1}$ and $W^{2}$ respectively.\\ \indent We introduced a few fundamental concepts in random dynamical system and nonlocal or fractional Laplacian operator next section.\end{linenomath*}
\section{Preliminaries}

\indent Consider $H_{1}=L^{2}(-1,1)$, which denotes the space for fast variables. While $H_{2}$ is a separable Hilbert space, which denotes the space for slow variables. Take $A_{\alpha}=-(-\Delta)^{\frac{\alpha}{2}}$, which represents the nonlocal fractional Laplacian operator. On the functions that are extended to $\mathbb{R}$ by zero just as in \cite{kwasnicki2012eigenvalues}, the above nonlocal fractional Laplacian is applied. Now we review a few basic notions in random systems.
\bd (\cite{arnold2013random})
\begin{linenomath*}Let $(\Omega, \mathcal{F},\mathbb{P})$ be a probability space and $\theta=\{\theta_{r}\}_{r\in\mathbb{R}}$ be a flow on $\Omega$, namely a mapping $$\theta:\mathbb{R}\times\Omega\rightarrow\Omega,$$ and satisfies the following conditions\\
 $\bullet$  $\theta_{0}=Id_{\Omega};$\\
   $\bullet$ $\theta_{r_{1}}\theta_{r_{2}}=\theta_{r_{1}+r_{2}},$ where $r_{1},r_{2}\in\mathbb{R};$\\
   $\bullet$ the mapping $(r,\omega)\mapsto \theta_{r}\omega$ is $(\mathcal{B(\mathbb{R})\otimes\mathcal{F},F})$-measurable, and $\theta_{r}\mathbb{P}=\mathbb{P}$ for all $r\in\mathbb{R}$. Then $\Theta=(\Omega, \mathcal{F},\mathbb{P},\theta)$ is known as a driving or metric dynamical system.\end{linenomath*}\ed
\begin{linenomath*}To achieve our objective, we take a special but very important driving dynamical system influenced by Wiener process. Let us consider $W(t)$ taking values in state space (Hilbert space) $\mathbb{H}$ be a two-sided Wiener process, which depends upon time $t$ and also varying with $t$. The defined Wiener process has zero value at $t=0$. Its sample path is the space of real continuous functions on $\mathbb{R}$, which is denoted by $C_{0} (\mathbb{R},\mathbb{H})$. Here we take the flow, which is also measurable and defined by $\theta=\{\theta_{r}\}_{r\in\mathbb{R}}$
$$\theta_{r}\omega=\omega(\cdot+r)-\omega(r),\indent \omega \in \Omega, \indent r\in\mathbb{R}.$$\indent On $\mathcal{B}(C_{0} (\mathbb{R},\mathbb{H}))$, the distribution of above process induces a probability measure, which is known as the Wiener measure. Instead of the whole $C_{0} (\mathbb{R},\mathbb{H})$, we consider an $\{\theta_{r}\}_{r\in\mathbb{R}}$-invariant subset $\Omega \subset C_{0} (\mathbb{R},\mathbb{H})$, whose probability-measure is one. With respect to $\Omega$ we also take the trace $\sigma$-algebra $\mathcal{F}$ of $\mathcal{B}(C_{0} (\mathbb{R},\mathbb{H}))$. Recall the concept of $\{\theta_{r}\}_{r\in\mathbb{R}}$-invariant, a set $\Omega$ is said to be $\{\theta_{r}\}_{r\in\mathbb{R}}$-invariant if $\theta_{r}\omega=\omega$ for $r\in\mathbb{R}$. Here we take the restriction of Weiner measure and we still represent it by $\mathbb{P}.$\\\indent On the state space $\mathbb{H}$ (which is often taken as a Hilbert space), the dynamics of a stochastic system over a flow $\theta$  is depicted by a cocycle in general.\end{linenomath*}
\bd (\cite{arnold2013random})
\begin{linenomath*}A cocycle $\phi$ is defined by a map:$$\phi:\mathbb{R}\times\Omega\times\mathbb{H}\rightarrow\mathbb{H},$$ which is $(\mathcal{B(\mathbb{R})\otimes\mathcal{F}\otimes\mathcal{B(\mathbb{H})},\mathcal{F}})$-measurable and satisfy the conditions \begin{align*}&\phi(0,\omega,y)=y,\\&\phi(r_{1}+r_{2},\omega,y)=\phi(r_{2},\theta_{r_{1}}\omega,\phi(r_{1},\omega,y)),\end{align*}for $y\in\mathbb{H}$, $\omega\in\Omega$ and $r_{1},r_{2}\in\mathbb{R}$. Then $\phi$, together with metric dynamical system $(\Omega,\mathcal{F},\mathbb{P},\theta)$, constructs a random dynamical system.\end{linenomath*}\ed
\begin{linenomath*}We say that a random dynamical system is differentiable (continuous) if $y\mapsto\phi(r,\omega,y)$ is differentiable (continuous) for $\omega\in\Omega$ and $r\geq0$. The metric space $(\mathbb{H},||\cdot||_{\mathbb{H}})$ contains a collection of closed and nonempty sets $M=\{M(\omega):\omega\in\Omega\}$, if for every $
y'\in\mathbb{H}$ the mapping defined below:$$\omega\mapsto
\mathop {\inf }\limits_{y \in M(\omega)}||y-y'||_{\mathbb{H}},$$ is a random variable. Then this family is said to be a random set.\end{linenomath*}
\bd
\begin{linenomath*}(\cite{duan2015introduction}) A random variable $y(\omega),$ which take values in $\mathbb{H}$, is known as stationary orbit (random fixed point) for a random system $\phi$ if$$\phi(r,\omega,y(\omega))=y(\theta_{r}\omega),\indent \indent a.s.$$for all $r$.\end{linenomath*}\ed
\begin{linenomath*}Now idea of random invariant manifold is given.\end{linenomath*}
\bd (\cite{fu2012slow})
\begin{linenomath*}If a random set $M=\{M(\omega):\omega\in \Omega\}$ satisfy the condition $$\phi(r,\omega,M(\omega))\subset M(\theta_{r}\omega),$$ for $\omega\in\Omega$ and $r\geq0$. Then $M$ is known as a random positively invariant set.\end{linenomath*}\ed
\begin{linenomath*}If $M=\{M(\omega):\omega\in \Omega\}$ can be represent by a Lipschitz mapping graph define as:$$h(\omega,\cdot):H_{2}\rightarrow H_{1},$$such that$$M(\omega)=\{(h(\omega,v),v):v\in H_{2}\},$$ then $M$ is said to be a Lipschitz random invariant manifold. Additionally, if for all $y\in \mathbb{H}$, there is an $y'\in M(\omega)$, such that $$||\phi(r,\omega,y)-\phi(r,\omega,y')||_{\mathbb{H}}\leq c_{1}(y,y',\omega)e^{c_{2}t}||y-y'||_{\mathbb{H}},\mbox{ for all }\omega\in \Omega,\mbox{ and }r\geq0.$$Where $c_{1}$ is a positive variable, which depend upon $y$ and $y'$, while  $-c_{2}$ is positive constant, then the random set $M(\omega)$ have exponential tracking property.\\
\indent Now we review the asymptotic behavior of eigenvalues of $A_{\alpha}$, $(\alpha\in (0,2))$ for the spectral problem $$(-\Delta)^{\frac{\alpha}{2}}\varphi(y)=\lambda\varphi(y),\indent y\in (-1,1),$$where $\varphi(\cdot)\in L^{2}(-1,1)$ is extended to $\mathbb{R}$ by zero just as in \cite{kwasnicki2012eigenvalues}. It is known that there exist an infinite sequence of eigenvalues $\lambda_{k}$ such that
$$0<\lambda_{1}<\lambda_{2}\leq \lambda_{3}\leq\cdot\cdot\cdot\leq\lambda_{k}\leq \cdot\cdot\cdot, \mbox{ for } k=1,2,3,\cdot\cdot\cdot.$$ and the corresponding eigenfunctions $\varphi_{k}$ form a complete orthonormal set in $L^{2}(-1,1)$.\end{linenomath*}\\
\textbf{Lemma 1.} \begin{linenomath*}(\cite{kwasnicki2012eigenvalues}) For the above spectral problem$$\lambda_{k}=(\frac{k\pi}{2}-\frac{(2-\alpha)\pi}{8})^{\alpha}+O(\frac{1}{k}),\indent (k\rightarrow\infty).$$\end{linenomath*}
\textbf{Lemma 2.}\begin{linenomath*} (\cite{yan2014approximation}) The Laplacian operator $A_{\alpha} $ is a sectorial operator, which satisfies the following upper-bound$$||e^{A_{\alpha}t}||_{L^{2}(-1,1)}\leq C e^{-\lambda_{1}t},\mbox{ for } t\geq0.$$Where $C$ is independent of $t$ and $\lambda_{1}$, also $C>0$ is constant.\\ \\
\indent Note that $(-A_{\alpha})^{-1}$ is a linear bounded operator on $L^{2}(-1,1)$, together with it is self-adjoint and compact operator.\end{linenomath*}
\section{Framework}

\begin{linenomath*} Framework includes a list of assumptions and conversion of stochastic dynamical system to random dynamical system.\\ Consider the system of stochastic evolutionary equations (1)-(2), which depend upon two time scales as described in previous section 1.
We assume the following hypothesis in the fast-slow system (1)-(2).\\\textbf{Hypothesis H1.} (Slow evolution): The operator $J$ generates $C_{0}$-semigroup $e^{Jt}$ on $H_{2}$, which satisfies $$||e^{Jt}v||_{2}\leq e^{\gamma_{2} t}||v||_{2},\indent t\leq 0,$$ with a constant $\gamma_{2}
>0$ and for all $v\in H_{2}.$\\\textbf{Hypothesis H2.} (Lipschitz condition): Nonlinearities $f$ and $g$ are continuous and differentiable functions ($C^{1}-smooth)$,$$f:\mathbb{H}\rightarrow H_{1},\indent g:\mathbb{H}\rightarrow H_{2},$$with$$f(0,0)=0=g(0,0), \indent \frac{\partial f}{\partial u^{\epsilon}}(0,0)=\frac{\partial f}{\partial v^{\epsilon}}(0,0)=0=\frac{\partial g}{\partial u^{\epsilon}}(0,0)=\frac{\partial g}{\partial v^{\epsilon}}(0,0).$$ For all $(u^{\epsilon}_{i},v^{\epsilon}_{i})^{T}\in H_{1}\times H_{2}$, $\exists$ a positive constant $K$, such that$$||f(u^{\epsilon}_{1},v^{\epsilon}_{1})-f(u^{\epsilon}_{2},v^{\epsilon}_{2})||\leq K (||u^{\epsilon}_{1}-u^{\epsilon}_{2}||_{1}+||v^{\epsilon}_{1}-v^{\epsilon}_{2}||_{2}),$$also$$||g(u^{\epsilon}_{1},v^{\epsilon}_{1})-g(u^{\epsilon}_{2},v^{\epsilon}_{2})||\leq K (||u^{\epsilon}_{1}-u^{\epsilon}_{2}||_{1}+||v^{\epsilon}_{1}-v^{\epsilon}_{2}||_{2}).$$\textbf{Hypothesis H3.} (Gap condition): In the system, the term $K$ (Lipschitz constant) of the functions $f$ and $g$ satisfies the bound  $$K<\frac{\lambda_{1}\gamma_{2}}{2\lambda_{1}+\gamma_{2}},$$ where $\lambda_{1}$ is the decay rate of Laplacian operator $A_{\alpha}$.\\ \indent Now we construct a random evolutionary system corresponding to the stochastic evolutionary system (1)-(2). For the achievement of this objective, we first set up the existence, after that verify the uniqueness of solution for the equations (1)-(2) and for Ornstein-Uhblenbeck equation, which is nonlocal equation.\\
\textbf{Lemma 3.} Let $W$ taking values in Hilbert space $\mathbb{H}$ be a two sided Weiner process (a time dependent Brownian motion) then under hypothesis $(H1-H3)$, the system (1)-(2) posses a unique mild solution.\\
\textbf{Proof.} We write the system (1)-(2) in the following form \begin{align}\label{sde} \left( {\begin{array}{*{20}{c}}
{\dot{u}^{\epsilon}}\\
\dot{v}^{\epsilon}
\end{array}} \right)=
\left( {\begin{array}{*{20}{c}}
{\frac{1}{\epsilon }{A_\alpha }}&0\\
0&J
\end{array}} \right)
\left( {\begin{array}{*{20}{c}}
{u^{\epsilon}}\\
v^{\epsilon}
\end{array}} \right)
+\left( {\begin{array}{*{20}{c}}
{\frac{1}{\epsilon}f(u^{\epsilon},v^{\epsilon})}\\
g(u^{\epsilon},v^{\epsilon})
\end{array}} \right)
+\left( {\begin{array}{*{20}{c}}
{\frac{\sigma_{1}}{\sqrt{\epsilon}}\dot{W}_{t}^{1}}\\
\sigma_{2}\dot{W}_{t}^{2}
\end{array}} \right)
\end{align}
By Lemma 1 and $\epsilon>0$, we know that for every $u$ in $D(\frac{1}{\epsilon}A_{\alpha})$, there exist $u^{*} \in \mathbb{H}^{*}$ (dual space of Hilbert space) such that $\langle\frac{1}{\epsilon}A_{\alpha}u,u^{*}\rangle <0$, hence $A_{\alpha}$ is dissipative and $\lambda I-A_{\alpha}$ is invertible for any $\lambda>0$. This implies the range of $\lambda I-A_{\alpha}$ is whole space $\mathbb{H}$.\\ So by Lumer-Phillips theorem \cite{pazy2012semigroups}, $\frac{1}{\epsilon}A_{\alpha}$ is infinitesimal generator of a semigroup which is strongly continuous. Considering assumption $(H1)$, we can say that the operator $$\left( {\begin{array}{*{20}{c}}
{\frac{1}{\epsilon }{A_\alpha }}&0\\
0&J
\end{array}} \right)$$ is also an infinitesimal generator of semigroup, which is strongly continuous. Since, $W$ taking values in Hilbert space $\mathbb{H}$ is a two sided Weiner process. Hence by using Theorem 7.2(\cite{da2014stochastic},p.188), the above system (1)-(2) posses a unique mild solution.  $\Box$\\ \textbf{Lemma 4.} The nonlocal stochastic equation $$ d\rho(t)=A_{\alpha}\rho(t)dt+\sigma_{1}dW_{t}^{1},\mbox{ in } H_{1},$$where $A_{\alpha}$ is the Laplacian operator, has the following solution$$ \rho(t)=e^{tA_{\alpha}}\rho_{0}+\sigma_{1}\int_{0}^{t}e^{A_{\alpha}(t-s)}dW_{s}^{1},\mbox{ for all } t\geq0 \mbox{ and } \alpha \in
(1,2),$$here $\rho(0)=\rho_{0}$ is $\mathcal{F}_{0}$-measurable.\\By an equation is in $H_{1}$ means that every term in the equation is in $H_{1}$.\\ \textbf{Remark 1.} \cite{bai2017slow} The above lemma is not hold for $\alpha \in (0,1]$ .\\ Now, let $\Psi_{1}=(\Omega_{1},\mathcal{F}_{1},\mathbb{P}_{1},\theta_{t}^{1})$ and $\Psi_{2}=(\Omega_{2},\mathcal{F}_{2},\mathbb{P}_{2},\theta_{t}^{2})$ are two driving dynamical systems as explained in section 2. Define$$\Psi=\Psi_{1}\times\Psi_{2}=(\Omega_{1}\times\Omega_{2},\mathcal{F}_{1}\otimes\mathcal{F}_{1},\mathbb{P}_{1}\times\mathbb{P}_{1},\theta_{t}^{1}\times\theta_{t}^{2})$$
and $$\theta_{t}\omega:=(\theta_{t}^{1}\omega_{1},\theta_{t}^{2}\omega_{2})^{T}, \mbox{ for } \omega:=(\omega_{1},\omega_{2})^{T}\in \Omega_{1}\times\Omega_{2}:=\Omega.$$ Consider the stochastic evolutionary equations
\begin{align}&d\eta^{\epsilon}(t)=
\frac{1}{\epsilon}A_{\alpha}\eta^{\epsilon}dt+\frac{\sigma_{1}}{\sqrt{\epsilon}}dW_{t}^{1},\\&d\xi(t)=J\xi dt+\sigma_{2}dW_{t}^{2}.\end{align}
 \textbf{Lemma 5.} (\cite{schmalfuss2008invariant})  Assume that the hypothesis (H1) holds. Then equations (5) and (6) have continuous stationary solutions $\eta^{\epsilon}(\theta_{t}^{1}\omega_{1})$ and $\xi(\theta_{t}^{2}\omega_{2})$ respectively.\\Introduce a random transformation \begin{align}\binom{U^{\epsilon}}{V^{\epsilon}}:=\nu(\omega,u^{\epsilon},v^{\epsilon})=\binom{u^{\epsilon}-\eta^{\epsilon}(\theta_{t}^{1}\omega_{1})}{v^{\epsilon}-\xi(\theta_{t}^{2}\omega_{2})},\end{align}so the original stochastic evolutionary system (1)-(2) is converted to a random evolutionary system as\\
\begin{align} &dU^{\epsilon}=\frac{1}{\epsilon}A_{\alpha}U^{\epsilon}dt+\frac{1}{\epsilon}f(U^{\epsilon}+\eta^{\epsilon}(\theta_{t}^{1}\omega_{1}), V^{\epsilon}+\xi(\theta_{t}^{2}\omega_{2}))dt,\\ &dV^{\epsilon}=JV^{\epsilon}dt+g(U^{\epsilon}+\eta^{\epsilon}(\theta_{t}^{1}\omega_{1}), V^{\epsilon}+\xi(\theta_{t}^{2}\omega_{2}))dt.\end{align}The nonlinear functions $f$ and $g$ in (8)-(9) satisfies the same Lipschitz condition. Let $Z^{\epsilon}(t,\omega,Z_{0})=(U^{\epsilon}(t,\omega,(U_{0},V_{0})^{T}),V^{\epsilon}(t,\omega,(U_{0},V_{0})^{T}))^{T}$ be the solution of (8)-(9) with initial value $Z_{0}:=(U^{\epsilon}(0),V^{\epsilon}(0))^{T}=(U_{0},V_{0})^{T}$. By the classical theory for evolutionary equations \cite{EvansL2W2015}, under the hypothesis (H1-H3) the system (8)-(9) with initial data has unique global solution for every $\omega=(\omega_{1},\omega_{2})^{T}\in\Omega=\Omega_{1}\times\Omega_{2}$ just like to \cite{chen2014slow}.\\Hence the solution mapping of random evolutionary system (8)-(9) $$\Phi^{\epsilon}(t,\omega,(U_{0},V_{0})^{T}):=(U^{\epsilon}(t,\omega,(U_{0},V_{0})^{T}),V^{\epsilon}(t,\omega,(U_{0},V_{0})^{T}))^{T},$$indicates a random dynamical system. Furthermore, $$\phi^{\epsilon}(t,\omega):=\Phi^{\epsilon}(t,\omega)+(\eta^{\epsilon}(\theta_{t}^{1}\omega_{1}),\xi(\theta_{t}^{2}\omega_{2})),\indent t\geq0,\indent \omega\in\Omega.$$is the random system, which is generated by the original stochastic system (1)-(2).\end{linenomath*}
\section{Slow manifold}
\begin{linenomath*}In this part, we construct the slow manifold existence and verify the exponential tracking property for the random evolutionary system (8)-(9). Define$$M^{\epsilon}(\omega)\triangleq\{(H^{\epsilon}(\omega,V),V)^{T}:V\in H_{2}\}.$$We prove that the random set $M^{\epsilon}(\omega)$ is invariant manifold defined as a graph of Lipschitz mapping with the help of Lyapunov-Perron method. To study system (8)-(9), we construct the following Banach spaces in term of functions with a weighted geometrically sup norm \cite{wanner1995linearization}. For any $\beta \in\mathbb{R}$:\\$\indent C_{\beta}^{1,-}=\{\chi:(-\infty,0]\rightarrow H_{1}$ is continuous and $\mathop {\mbox{sup} }\limits_{t\leq 0}e^{-\beta t}||\chi(t)||_{1}<\infty\},$\\$\indent C_{\beta}^{2,-}=\{\chi:(-\infty,0]\rightarrow H_{2}$ is continuous and $\mathop {\mbox{sup} }\limits_{t\leq 0}e^{-\beta t}||\chi(t)||_{2}<\infty\},$\\ with the norms\\$\indent ||\chi||_{C_{\beta}^{1,-}}=\mathop {\mbox{sup} }\limits_{t\leq 0}e^{-\beta t}||\chi(t)||_{1},$ and $||\chi||_{C_{\beta}^{2,-}}=\mathop {\mbox{sup} }\limits_{t\leq 0}e^{-\beta t}||\chi(t)||_{2}.$\\Similarly we define,\\
$\indent C_{\beta}^{1,+}=\{\chi:[0,\infty)\rightarrow H_{1}$ is continuous and $\mathop {\mbox{sup} }\limits_{t\geq 0}e^{-\beta t}||\chi(t)||_{1}<\infty\},$\\$\indent C_{\beta}^{2,+}=\{\chi:[0,\infty)\rightarrow H_{2}$ is continuous and $\mathop {\mbox{sup} }\limits_{t\geq 0}e^{-\beta t}||\chi(t)||_{2}<\infty\},$\\ with the norms\\$\indent ||\chi||_{C_{\beta}^{1,+}}=\mathop {\mbox{sup} }\limits_{t\geq 0}e^{-\beta t}||\chi(t)||_{1},$ and $||\chi||_{C_{\beta}^{2,+}}=\mathop {\mbox{sup} }\limits_{t\geq 0}e^{-\beta t}||\chi(t)||_{2}.$\\ Let $C_{\beta}^{\mp}$ be the product Banach spaces $C_{\beta}^{\mp}:=C_{\beta}^{1,\mp}\times C_{\beta}^{2,\mp}$, with the norm$$||z||_{C_{\beta}^{\mp}}=||u^{\epsilon}||_{C_{\beta}^{1,\mp}}+||v^{\epsilon}||_{C_{\beta}^{2,\mp}},\indent z=(u^{\epsilon},v^{\epsilon})^{T}\in C_{\beta}^{\mp}.$$For convenience, let $\mu$ be a number, which satisfy\begin{align}\label{sde} 0<\mu=\frac{\gamma_{2}}{2\lambda_{1}+\gamma_{2}}<1,\mbox{ and } K<\mu\lambda_{1}<\lambda_{1}, \mbox{ also } -\mu+\lambda_{1}>K.\end{align}Now, for the achievement of our goal, the following Lemma from \cite{chen2014slow} is needed.\\\textbf{Lemma 6.} Consider the so-called Lyapunov-Perron equation\\
$H^{\epsilon}(\omega,V_{0})=\frac{1}{\epsilon}\int_{-\infty}^{0}e^{-\frac{A_{\alpha}s}{\epsilon}}f(U^{\epsilon}(s,\omega;V_{0})+\eta^{\epsilon}(\theta_{s}^{1}\omega_{1}),V^{\epsilon}(s,\omega;V_{0})+\xi(\theta_{s}^{2}\omega_{2}))ds,$\\
for $V_{0}\in H_{2}$, where $U^{\epsilon}(s,\omega;V_{0})$ and $V^{\epsilon}(s,\omega;V_{0})$ are the solutions of system (8)-(9) with initial value $Z_{0}=(U_{0},V_{0})$ satisfying the form \begin{align*}&\binom{U^{\epsilon}(t,\omega;V_{0})}{V^{\epsilon}(t,\omega;V_{0})}\\&=\binom{\frac{1}{\epsilon}\int_{-\infty}^{t}e^{A_{\alpha}(t-s)/\epsilon}f(U^{\epsilon}(s,\omega;V_{0})+\eta^{\epsilon}(\theta_{s}^{1}\omega_{1}),V^{\epsilon}(s,\omega;V_{0})+\xi(\theta_{s}^{2}\omega_{2}))ds}{e^{Jt}V_{0}+\int_{0}^{t}e^{J(t-s)}g(U^{\epsilon}(s,\omega;V_{0})+\eta^{\epsilon}(\theta_{s}^{1}\omega_{1}),V^{\epsilon}(s,\omega;V_{0})+\xi(\theta_{s}^{2}\omega_{2}))ds}.\end{align*}
Then $Z^{\epsilon}(.,\omega)=(U^{\epsilon}(.,\omega),V^{\epsilon}(.,\omega))^{T}\in C_{\beta}^{-}.$\\
\textbf{Theorem 1.} (Slow manifold)
Suppose that (H1-H3) hold and $\epsilon>0$ is small enough. Then random evolutionary system (8)-(9) possessing a Lipschitz random slow manifold $M^{\epsilon}(\omega)$ represented by a graph$$M^{\epsilon}(\omega)=\{(H^{\epsilon}(\omega,V),V)^{T}:V\in H_{2}),$$where $$H^{\epsilon}(\cdot,\cdot):\Omega\times H_{2}\rightarrow H_{1},$$ is a graph mapping with Lipschitz constant satisfying$$Lip H^{\epsilon}(\omega,\cdot)\leq\frac{K}{(-\mu+\lambda_{1})[1-K(\frac{1}{-\mu+\lambda_{1}}+\frac{\epsilon}{\epsilon\gamma_{2}+\mu})]}.$$
\textbf{Remark 2.} Since there is a relation among the solutions of stochastic evolutionary system (1)-(2) and random system (8)-(9), so if the system (1)-(2) satisfies the conditions described in Theorem 1, then it also possess Lipschitz random invariant manifold$$\breve{M}^{\epsilon}(\omega)=M^{\epsilon}(\omega)+(\eta^{\epsilon}(\omega_{1}),\xi(\omega_{2}))=\{(\breve{H}^{\epsilon}(\omega,V),V)^{T}:V\in H_{2}\},$$where$$\breve{H}^{\epsilon}(\omega,V)=H^{\epsilon}(\omega,V)+(\eta^{\epsilon}(\omega_{1}),0)).$$
\textbf{Remark 3.} If $f$ and $g$ are $C^{k}$ and a large spectrum gap condition holds then it can be proved similar to \cite{duan2004smooth}, that manifolds are $C^{k}$-smooth.\\
\textbf{Proof.} To establish an invariant manifold corresponding to the random dynamical system (8)-(9), we first consider the integral equation
\begin{align}&\binom{U^{\epsilon}(t)}{V^{\epsilon}(t)}=\binom{\frac{1}{\epsilon}\int_{-\infty}^{t}e^{A_{\alpha}(t-s)/\epsilon}f(U^{\epsilon}(s)+\eta^{\epsilon}(\theta_{s}^{1}\omega_{1}),V^{\epsilon}(s)+\xi(\theta_{s}^{2}\omega_{2}))ds}{e^{Jt}V_{0}+\int_{0}^{t}e^{J(t-s)}g(U^{\epsilon}(s)+\eta^{\epsilon}(\theta_{s}^{1}\omega_{1}),V^{\epsilon}(s)+\xi(\theta_{s}^{2}\omega_{2}))ds},t\leq0.\end{align} \textbf{\indent Step 1.} Let $Z^{\epsilon}(t,\omega,Z_{0})=(U^{\epsilon}(t,\omega,V_{0}) ,V^{\epsilon}(t,\omega,V_{0}))^{T}$ be the solution of (11) with initial value $Z_{0}=(U_{0},V_{0})$. We prove that $Z^{\epsilon}(t,\omega,Z_{0})$ is the unique solution of (11), by using of Banach fixed point theorem. For this,\\define two operators $\mathfrak{I}_{1}^{\epsilon}:C_{\beta}^{-}\rightarrow C_{\beta}^{1,-}$ and $\mathfrak{I}_{2}^{\epsilon}:C_{\beta}^{-}\rightarrow C_{\beta}^{2,-}$ by means of $$\mathfrak{I}_{1}^{\epsilon}(z(\cdot))[t]=\frac{1}{\epsilon}\int_{-\infty}^{t}e^{A_{\alpha}(t-s)/\epsilon}f(u(s)+\eta^{\epsilon}(\theta_{s}^{1}\omega_{1}),v(s)+\xi(\theta_{s}^{2}\omega_{2}))ds,$$
$$\mathfrak{I}_{2}^{\epsilon}(z(\cdot))[t]=e^{Jt}V_{0}+\int_{0}^{t}e^{J(t-s)}g(u(s)+\eta^{\epsilon}(\theta_{s}^{1}\omega_{1}),v(s)+\xi(\theta_{s}^{2}\omega_{2}))ds,$$ for $t\leq0$, and the Lyapunov-Perron transform $\mathfrak{I}^{\epsilon}$ is$$\mathfrak{I}^{\epsilon}(z(\cdot))=\binom{\mathfrak{I}_{1}^{\epsilon}(z(\cdot))}{\mathfrak{I}_{2}^{\epsilon}(z(\cdot))}^{T}=(\mathfrak{I}_{1}^{\epsilon}(z(\cdot)),\mathfrak{I}_{2}^{\epsilon}(z(\cdot))).$$
For verify that $\mathfrak{I}^{\epsilon}$ maps $C_{\beta}^{-}$ into itself. Take $z=(u,v)^{T}\in C_{\beta}^{-}$, such that:
\begin{align*}||\mathfrak{I}_{1}^{\epsilon}(z(\cdot))[t]||_{C_{\beta}^{1,-}}&=||\frac{1}{\epsilon}\int_{-\infty}^{t}e^{A_{\alpha}(t-s)/\epsilon}f(u(s)+\eta^{\epsilon}(\theta_{s}^{1}\omega_{1}),v(s)+\xi(\theta_{s}^{2}\omega_{2}))ds||_{1}\\
&\leq\frac{1}{\epsilon}\mathop { \mbox{sup} }\limits_{t\leq 0}\{e^{-\beta (t-s)}\int_{-\infty}^{t}e^{-\lambda_{1}(t-s)/\epsilon}||f(u(s)+\eta^{\epsilon}(\theta_{s}^{1}\omega_{1}),v(s)\\&\indent+\xi(\theta_{s}^{2}\omega_{2}))||_{1}ds\}\\
&\leq\frac{K}{\epsilon}\mathop {\mbox{sup} }\limits_{t\leq 0}\{e^{-\beta (t-s)}\int_{-\infty}^{t}e^{-\lambda_{1}(t-s)/\epsilon}(||u(s)||_{1}+||v(s)||_{2})ds\}+ C_{1}\\ &\leq\frac{K}{\epsilon}\mathop {\mbox{sup} }\limits_{t\leq 0}\{\int_{-\infty}^{t}e^{(-\lambda_{1}/\epsilon-\beta)(t-s)}ds\}||z||_{C_{\beta}^{-}}+ C_{1}\\&=\frac{K}{\epsilon\beta+\lambda_{1}}||z||_{C_{\beta}^{-}}+ C_{1},\end{align*}
and
\begin{align*}||\mathfrak{I}_{2}^{\epsilon}(z(\cdot))[t]||_{C_{\beta}^{2,-}}&=||e^{Jt}V_{0}+\int_{0}^{t}e^{J(t-s)}g(u(s)+\eta^{\epsilon}(\theta_{s}^{1}\omega_{1}),v(s)+\xi(\theta_{s}^{2}\omega_{2}))ds||_{2}\\
&\leq \mathop {\mbox{sup} }\limits_{t\leq 0}\{e^{-\beta (t-s)}\int_{t}^{0}e^{\gamma_{2}(t-s)}||g(u(s)+\eta^{\epsilon}(\theta_{s}^{1}\omega_{1}),v(s)\\&\indent+\xi(\theta_{s}^{2}\omega_{2}))||_{2}ds\}+\mathop {\mbox{sup} }\limits_{t\leq 0}\{e^{-\beta t}e^{\gamma_{2} t}||V_{0}||_{2}\}\\
&\leq K \mathop {\mbox{sup} }\limits_{t\leq 0}(\int_{t}^{0}e^{(\gamma_{2}-\beta)(t-s)}ds)(||u(s)||_{1}+||v(s)||_{2})+C_{2}+||V_{0}||_{2}\\&=\frac{ K}{\gamma_{2}-\beta}||z||_{C_{\beta}^{-}}+C_{2}+||V_{0}||_{2}\\&=\frac{K}{\gamma_{2} -\beta}||z||_{C_{\beta}^{-}}+C_{3}.\end{align*} Hence, by combining with the definition of $\mathfrak{I}$, we obtain that\\$$||\mathfrak{I}^{\epsilon}(z)||_{C_{\beta}^{-}}\leq[\frac{K}{\epsilon\beta+\lambda_{1}}+\frac{ K}{\gamma_{2}-\beta}]||z||_{C_{\beta}^{-}}+C_{4},$$ where $C_{i}, i=1,2,3,4$ are constants and $$\upsilon(K,\beta,\lambda_{1},\gamma_{2},\epsilon)= \frac{K}{\epsilon\beta+\lambda_{1}}+\frac{K}{\gamma_{2}-\beta}.$$ This implies that $\mathfrak{I}^{\epsilon}$ maps $C_{\beta}^{-}$ into itself.\\Furthermore, to show that the map $\mathfrak{I}^{\epsilon}$ is contractive.\\Take $z=(u,v)^{T},\bar{z}=(\bar{u},\bar{v})^{T}\in C_{\beta}^{-}$, \begin{align*}||\mathfrak{I}_{1}^{\epsilon}(z)-\mathfrak{I}_{1}^{\epsilon}(\bar{z})||_{C_{\beta}^{1,-}}&\leq\frac{1}{\epsilon}\mathop {\mbox{sup} }\limits_{t\leq 0}\{e^{-\beta (t-s)}\int_{-\infty}^{t}e^{-\lambda_{1}(t-s)/\epsilon}||f(u(s)+\eta^{\epsilon}(\theta_{s}^{1}\omega_{1}),v(s)\\& \indent +\xi(\theta_{s}^{2}\omega_{2}))-f(\bar{u}(s)+\eta^{\epsilon}(\theta_{s}^{1}\omega_{1}),\bar{v}(s)+\xi(\theta_{s}^{2}\omega_{2}))||_{1}ds\}\\
&\leq\frac{K}{\epsilon}\mathop {\mbox{sup} }\limits_{t\leq 0}\{e^{-\beta (t-s) }\int_{-\infty}^{t}e^{-\lambda_{1}(t-s)/\epsilon}(||u(s)-\bar{u}(s)||_{1}\\&\indent+||v(s)-\bar{v}(s)||_{2})ds\}\\
&\leq\frac{K}{\epsilon}\mathop {\mbox{sup} }\limits_{t\leq 0}\{\int_{-\infty}^{t}e^{(\frac{-\lambda_{1}}{\epsilon}-\beta)(t-s)}ds\}||z-\bar{z}||_{C_{\beta}^{-}}\\
&=\frac{K}{\epsilon\beta+\lambda_{1}}||z-\bar{z}||_{C_{\beta}^{-}}.\end{align*} Similarly
\begin{align*}||\mathfrak{I}_{2}^{\epsilon}(z)-\mathfrak{I}_{2}^{\epsilon}(\bar{z})||_{C_{\beta}^{2,-}}&\leq K \mathop {\mbox{sup} }\limits_{t\leq 0}\{\int_{t}^{0}e^{\gamma_{2}(t-s)}e^{-\beta(t-s)}ds\}||z-\bar{z}||_{C_{\beta}^{-}}\\
&\leq K \mathop {\mbox{sup} }\limits_{t\leq 0}\{\int_{t}^{0}e^{(\gamma_{2}-\beta)(t-s)}ds\}||z-\bar{z}||_{C_{\beta}^{-}}\\&=\frac{ K}{\gamma_{2}-\beta}||z-\bar{z}||_{C_{\beta}^{-}},\end{align*} which implies that\begin{align*} &||\mathfrak{I}^{\epsilon}(z)-\mathfrak{I}^{\epsilon}(\bar{z})||_{C_{\beta}^{-}}\leq\upsilon(K,\beta,\lambda_{1},\gamma_{2},\epsilon)||z-\bar{z}||_{C_{\beta}^{-}},\end{align*}
where\begin{align*}&\upsilon(K,\beta,\lambda_{1},\gamma_{2},\epsilon)=\frac{K}{\epsilon\beta+\lambda_{1}}+\frac{K}{\gamma_{2}-\beta},\end{align*}
Notice the hypothesis (H3) holds, take $\beta=-\frac{\mu}{\epsilon}$, and then
$$\upsilon(K,\beta,\lambda_{1},\gamma_{2},\epsilon)\rightarrow\frac{K}{-\mu+\lambda_{1}} \mbox{ as } \epsilon\rightarrow0.$$
Therefore, there is a enough small parameter $\epsilon_{0}\rightarrow 0$ such that$$0<\upsilon(K,\beta,\lambda_{1},\gamma_{2},\epsilon)< 1,\mbox{ for }\epsilon \in (0,\epsilon_{0}).$$ Hence, the mapping $\mathfrak{I}^{\epsilon}$ is contractive in $C_{-\frac{\mu}{\epsilon}}^{-}$, and then (11) posses a unique solution$$Z^{\epsilon}(t,\omega,V_{0})=(U^{\epsilon}(t,\omega,V_{0}),V^{\epsilon}(t,\omega,V_{0}))^{T} \mbox{ in } C_{-\frac{\mu}{\epsilon}}^{-}.$$Furthermore,\begin{align*}||Z^{\epsilon}(\cdot,\omega,V_{1})-Z^{\epsilon}(\cdot,\omega,V_{2})||_{C_{-\frac{\mu}{\epsilon}}^{-}} &=||U^{\epsilon}(\cdot,\omega,V_{1})-U^{\epsilon}(\cdot,\omega,V_{2})||_{C_{-\frac{\mu}{\epsilon}}^{1,-}}+||V^{\epsilon}(\cdot,\omega,V_{1})\\&\indent-V^{\epsilon}(\cdot,\omega,V_{2})||_{C_{-\frac{\mu}{\epsilon}}^{2,-}}\\
&\leq\frac{K}{\epsilon\beta+\lambda_{1}}||Z^{\epsilon}(\cdot,\omega,V_{1})-Z^{\epsilon}(\cdot,\omega,V_{2})||_{C_{-\frac{\mu}{\epsilon}}^{-}}+\frac{ K}{\gamma_{2}-\beta}\\&\indent\times||Z^{\epsilon}(\cdot,\omega,V_{1})-Z^{\epsilon}(\cdot,\omega,V_{2})||_{C_{-\frac{\mu}{\epsilon}}^{-}} +||V_{1}-V_{2}||_{2}\\
&=\upsilon(K,\beta,\lambda_{1},\gamma_{2},\epsilon)||Z^{\epsilon}(\cdot,\omega,V_{1})-Z^{\epsilon}(\cdot,\omega,V_{2})||_{C_{-\frac{\mu}{\epsilon}}^{-}}\\&\indent+||V_{1}-V_{2}||_{2}.\end{align*}
Thus, we obtain the following upper-bound\begin{align}||Z^{\epsilon}(\cdot,\omega,V_{1})-Z^{\epsilon}(\cdot,\omega,V_{2})||_{C_{-\frac{\mu}{\epsilon}}^{-}}\leq\frac{1}{1-\upsilon(K,\beta,\lambda_{1},\gamma_{2},\epsilon)}||V_{1}-V_{2}||_{2},\end{align} for all $\omega \in \Omega, V_{1}, V_{2} \in H_{2}.$\\\textbf{\indent Step 2.} With the help of unique solution, which described in (Step 1), we construct the mapping $H^{\epsilon}.$\begin{equation}\label{sde}H^{\epsilon}(\omega, V_{0})=\frac{1}{\epsilon}\int_{-\infty}^{0}e^{-A_{\alpha}s/\epsilon}f(U^{\epsilon}(s,\omega,V_{0})+\eta^{\epsilon}(\theta_{s}^{1}\omega_{1}),V^{\epsilon}(s,\omega,V_{0})+\xi(\theta_{s}^{2}\omega_{2}))ds,\end{equation}
then owing to (12), below upper-bound is obtained.
$$||H^{\epsilon}(\omega,V_{1})-H^{\epsilon}(\omega,V_{2})||_{1}\leq\frac{K}{\epsilon\beta+\lambda_{1}}\frac{1}{[1-\upsilon(K,\beta,\lambda_{1},\gamma_{2},\epsilon)]}||V_{1}-V_{2}||_{2},$$
for all $V_{1},V_{2} \in H_{2}$ and $\omega \in \Omega$. This implies that$$||H^{\epsilon}(\omega,V_{1})-H^{\epsilon}(\omega,V_{2})||_{1}\leq\frac{K}{-\mu+\lambda_{1}}\frac{1}{[1-\upsilon(K,\beta,\lambda_{1},\gamma_{2},\epsilon)]}||V_{1}-V_{2}||_{2},$$
for all $V_{1},V_{2} \in H_{2}$ and $\omega \in \Omega$. Then by (Lemma 5) it follows$$M^{\epsilon}(\omega)=\{(H^{\epsilon}(\omega,V),V)^{T}:V\in H_{2}\}.$$ \textbf{\indent Step 3.} In this step, for proving $M^{\epsilon}(\omega)$ is random set, we prove that for any $z=(u,v)^{T}\in \mathbb{H}=H_{1}\times H_{2}$,\begin{equation}\label{sde}\omega\mapsto \mathop {\mbox{inf} }\limits_{z'\in \mathbb{H}}||(u,v)^{T}-(H^{\epsilon}(\omega,\mathfrak{I}z'),\mathfrak{I}z')^{T}||,\end{equation}is measurable (Theorem III.9 in Casting and Valadier \cite{castaing1977convex}, p. 67). Let separable space $\mathbb{H}$ has a countable dense set, say, $\mathbb{H}_{c}$. Then right hand side of (14) get equality to \begin{equation}\label{sde}\mathop {\mbox{inf} }\limits_{z'\in \mathbb{H}_{c}}||(u,v)^{T}-(H^{\epsilon}(\omega,\mathfrak{I}z'),\mathfrak{I}z')^{T}||,\end{equation} which immediately follows from the continuity of $H^{\epsilon}(\omega,\cdot)$. Under infimum of (14) the measurability of any expression follows since $\omega\mapsto H^{\epsilon}(\omega,\mathfrak{I}z')$ is measurable for every $z' \in \mathbb{H}.$\\ \textbf{\indent Step 4.} Now it only remains to show that the random set $M^{\epsilon}(\omega)$ is positively invariant. That is, we prove that for each $Z_{0} =(U_{0},V_{0})^{T}\in M^{\epsilon}(\omega), Z^{\epsilon}(s,\omega,Z_{0}) \in M^{\epsilon}(\theta_{s}\omega)$, for all $s\geq 0.$ Finally for every fixed $s\geq0,$ we may observe that $Z^{\epsilon}(t+s,\omega,Z_{0})$ is a solution of \begin{align*}
&dU^{\epsilon}=\frac{1}{\epsilon}A_{\alpha}U^{\epsilon}dt+\frac{1}{\epsilon}f(U^{\epsilon}+\eta^{\epsilon}(\theta_{t}^{1}\omega_{1}), V^{\epsilon}+\xi(\theta_{t}^{2}\omega_{2}))dt,\\&dV^{\epsilon}=JV^{\epsilon}dt+g(U^{\epsilon}+\eta^{\epsilon}(\theta_{t}^{2}\omega_{2}), V^{\epsilon}+\xi(\theta_{t}^{2}\omega_{2}))dt,\end{align*} with initial data $Z(0)=(U(0),V(0))^{T}=Z^{\epsilon}(s,\omega,Z_{0})$.\\ Thus, $Z^{\epsilon}(t+s,\omega,Z_{0})=Z^{\epsilon}(t,\theta_{s}\omega,Z^{\epsilon}(s,\omega,Z_{0}))$. Since $Z^{\epsilon}(\cdot,\omega,Z_{0})$ in $C_{-\frac{\mu}{\epsilon}}^{-}$, then $Z^{\epsilon}(t,\theta_{s}\omega,Z^{\epsilon}(s,\omega,Z_{0}))\in C_{-\frac{\mu}{\epsilon}}^{-}$. So, $Z^{\epsilon}(s,\omega,Z_{0}) \in M^{\epsilon}(\theta_{s}\omega)$. This completes the proof.   $\Box$\\
\textbf{Theorem 2.} (Exponential tracking property) Suppose that assumptions (H1-H3) hold. Then Lipschitz invariant manifold for random stochastic evolutionary system (8)-(9) as established in Theorem 1, for $\epsilon>0$ small enough, posses the exponential tracking property in the sense:
For every $Z_{0}=(U_{0},V_{0}) \in \mathbb{H}$, there exist $\tilde{Z}_{0}=(\tilde{U}_{0},\tilde{V}_{0})\in M^{\epsilon}(\omega)$ such that
$$||\Phi^{\epsilon}(t,\omega,Z_{0})-\Phi^{\epsilon}(t,\omega,\tilde{Z}_{0})||\leq\mathfrak{C}_{m}e^{-\mathfrak{C}_{n}t}||Z_{0}-\tilde{Z}_{0}||,\indent t\geq0,$$
here $\mathfrak{C}_{m}$ and $\mathfrak{C}_{n}$ are positive constants.\\
\textbf{Remark 4.} By the relationship among the solutions of (1)-(2) and (8)-(9), if system (8)-(9) has an exponential tracking manifold so has the original stochastic system (1)-(2).\\
\textbf{Remark 5.} For any solution of $Z^{\epsilon}=(U^{\epsilon},V^{\epsilon})$ of (8)-(9), there is an orbit $\bar{Z}^{\epsilon}=(\bar{U}^{\epsilon},\bar{V}^{\epsilon})$ on $M^{\epsilon}$ satisfying the equation
$$\dot{\bar{V}}^{\epsilon}=J\bar{V}^{\epsilon}+g(H^{\epsilon}(\theta_{t}\omega,\bar{V}^{\epsilon})+\eta^{\epsilon}(\theta_{t}^{1}\omega_{1}),\bar{V}^{\epsilon}+\xi(\theta_{t}^{2}\omega_{2}))+\sigma_{2}W_{t}^{2}$$
such that $$||Z^{\epsilon}-\bar{Z}^{\epsilon}||\leq\frac{e^{-\frac{\mu}{\epsilon}t}}{1-K(\frac{1}{-\mu+\lambda_{1}}+\frac{\epsilon}{\epsilon\gamma_{2}+\mu})}
||Z_{0}-\bar{Z}_{0}||, \indent t\geq0,$$where$$Z_{0}=(U^{\epsilon}(0),V^{\epsilon}(0)),\indent \bar{Z}_{0}=(\bar{U}^{\epsilon}(0),\bar{V}^{\epsilon}(0)).$$
\textbf{Remark 6.} A simple system may not posses the exponential tracking property, because absence of global Lipschitz condition. For example, the system $\dot{u} = \frac{1}{\epsilon}(-�u+2 \sin^{2} u)+\cdot\cdot\cdot$ and $\dot{v} = 0+v^{2} +\cdot\cdot\cdot$, has a slow manifold, namely $u = 0$. It is attractive, but only in the finite domain $u< 0.55457$ , not attractive for all $H$ because absence of global Lipschitz condition.\\
\textbf{Proof.} Suppose that there are two solutions for random system (8)-(9).\\ Say, $Z^{\epsilon}(t)=(U^{\epsilon}(t),V^{\epsilon}(t))^{T}$ and $\tilde{Z}^{\epsilon}(t)=(\tilde{U}^{\epsilon}(t),\tilde{V}^{\epsilon}(t))^{T}$. \\Then $\mathfrak{Z}^{\epsilon}(t)=\tilde{Z}^{\epsilon}(t)-Z^{\epsilon}(t):=(X^{\epsilon}(t),Y^{\epsilon}(t))^{T}$ satisfies the equations
\begin{align}&dU^{\epsilon}=\frac{1}{\epsilon}A_{\alpha}U^{\epsilon}dt+\frac{1}{\epsilon}F(X^{\epsilon},Y^{\epsilon},\eta^{\epsilon}(\theta_{t}^{1}\omega_{1}),\xi(\theta_{t}^{2}\omega_{2}))dt,\\
&dV^{\epsilon}=JV^{\epsilon}dt+G(X^{\epsilon},Y^{\epsilon},\eta^{\epsilon}(\theta_{t}^{1}\omega_{1}),\xi(\theta_{t}^{2}\omega_{2}))dt,\end{align}
where\\
$F(X^{\epsilon},Y^{\epsilon},\eta^{\epsilon}(\theta_{t}^{1}\omega_{1}),\xi(\theta_{t}^{2}\omega_{2}))$
\indent \indent \indent$$=f(X^{\epsilon}+U^{\epsilon}+\eta^{\epsilon}(\theta_{t}^{1}\omega_{1}),Y^{\epsilon}+V^{\epsilon}+\xi(\theta_{t}^{2}\omega_{2}))-f(U^{\epsilon}+\eta^{\epsilon}(\theta_{t}^{1}\omega_{1}),V^{\epsilon}+\xi(\theta_{t}^{2}\omega_{2})),$$
and\\
$G(X^{\epsilon},Y^{\epsilon},\eta^{\epsilon}(\theta_{t}^{1}\omega_{1}),\xi(\theta_{t}^{2}\omega_{2}))$
\indent \indent \indent$$=g(X^{\epsilon}+U^{\epsilon}+\eta^{\epsilon}(\theta_{t}^{1}\omega_{1}),Y^{\epsilon}+V^{\epsilon}+\xi(\theta_{t}^{2}\omega_{2}))-g(U^{\epsilon}+\eta^{\epsilon}(\theta_{t}^{1}\omega_{1}),V^{\epsilon}+\xi(\theta_{t}^{2}\omega_{2})).$$
First of all, we claim that $\mathfrak{Z}^{\epsilon}(t)=((X^{\epsilon}(t),Y^{\epsilon}(t))^{T}$ is a solution of (16)-(17) in $C_{\beta}^{+}$ for $\beta=-\frac{\mu}{\epsilon}$
if
\begin{equation}\label{sde}\binom{X^{\epsilon}(t)}{Y^{\epsilon}(t)}=\binom{e^{A_{\alpha}t/\epsilon}X^{\epsilon}(0)+\frac{1}{\epsilon}\int_{0}^{t}e^{A_{\alpha}(t-s)/\epsilon}F(X^{\epsilon}(s),Y^{\epsilon}(s),\eta^{\frac{1}{\epsilon}}(\theta_{s}^{1}\omega_{1}),\xi(\theta_{s}^{2}\omega_{2}))ds}{\int_{+\infty}^{t}e^{J(t-s)}G(X^{\epsilon}(s),Y^{\epsilon}(s),\eta^{\frac{1}{\epsilon}}(\theta_{s}^{1}\omega_{1}),\xi(\theta_{s}^{2}\omega_{2}))ds}.\end{equation}
Similar to (Theorem 1), above claim is verified by using variation of constants formula. Since the verification has same steps as in (Theorem 1), so we omit here. Furthermore, to prove that (18) has solutions $(X^{\epsilon},Y^{\epsilon})^{T}$ in $C_{\beta}^{+}$ with $(X^{\epsilon}(0),Y^{\epsilon}(0))^{T}=(X(0),Y(0))^{T}$ such that $$(\tilde{U}_{0},\tilde{V}_{0})^{T}=(X(0),Y(0))^{T}+(U_{0},V_{0})^{T}\in M^{\epsilon}(\omega).$$Recall that (as it is proved in Theorem 1)
$$(\tilde{U}_{0},\tilde{V}_{0})^{T}\in M^{\epsilon}(\omega)$$ if and only if $$\tilde{U}_{0}=\frac{1}{\epsilon}
\int_{ - \infty }^{0} e^{A_{\alpha}(-s)}F(X^{\epsilon}(s,\tilde{V}_{0}),Y^{\epsilon}(s,\tilde{V}_{0}),\eta^{\epsilon}(\theta_{s}^{1}\omega_{1}),\xi(\theta_{s}^{2}\omega_{2}))ds.$$
Since $$(\tilde{U}_{0},\tilde{V}_{0})^{T}=(X(0),Y(0))^{T}+(U_{0},V_{0})^{T}.$$ So it follows that
$$(\tilde{U}_{0},\tilde{V}_{0})^{T}=(X(0),Y(0))^{T}+(U_{0},V_{0})^{T}\in M^{\epsilon}(\omega)$$ if and only if \\ \begin{align*}X(0)+U_{0}=&\frac{1}{\epsilon}\int_{-\infty}^{0}e^{A_{\alpha}(-s)}F(X^{\epsilon}(s,Y(0)\\&+{V}_{0}),Y^{\epsilon}(s,Y(0)+{V}_{0}),\eta^{\epsilon}(\theta_{s}^{1}\omega_{1}),\xi(\theta_{s}^{2}\omega_{2}))ds
\\=&H^{\epsilon}(\omega,Y(0)+V_{0}).\end{align*} In other words $$(\tilde{U}_{0},\tilde{V}_{0})^{T}=(X(0),Y(0))^{T}+(U_{0},V_{0})^{T}\in M^{\epsilon}(\omega)$$ if and only if\begin{align}X(0)=-U_{0}+H^{\epsilon}(\omega,Y(0)+V_{0}).\end{align} For every $\mathfrak{Z}=(X,Y)^{T}\in C_{\beta}^{+}$ for $\beta=-\frac{\mu}{\epsilon}$ and $t\geq0$ define
\begin{align*} &\mathfrak{K}_{1}^{\epsilon}(\mathfrak{Z}(\cdot))[t]:=e^{A_{\alpha}t/\epsilon}X(0)+\frac{1}{\epsilon}\int_{0}^{t}e^{A_{\alpha}(t-s)/\epsilon}F(X(s),Y(s),\eta^{\epsilon}(\theta_{s}^{1}\omega_{1}),\xi(\theta_{s}^{2}\omega_{2}))ds,\\
&\mathfrak{K}_{2}^{\epsilon}(\mathfrak{Z}(\cdot))[t]:=\int_{+\infty}^{t}e^{J(t-s)}G(X(s),Y(s),\eta^{\epsilon}(\theta_{s}^{1}\omega_{1}),\xi(\theta_{s}^{2}\omega_{2}))ds,\end{align*}
where $X(0)$ is given in (19).\\Introduce Lyapunov-Perron transform $\mathfrak{K}:C_{-\frac{\mu}{\epsilon}}^{+}\rightarrow C_{-\frac{\mu}{\epsilon}}^{+}$ as
$$\mathfrak{K}^{\epsilon}(\mathfrak{Z}(\cdot))=\binom{\mathfrak{K}_{1}^{\epsilon}(\mathfrak{Z}(\cdot))}{\mathfrak{K}_{2}^{\epsilon}(\mathfrak{Z}(\cdot))}^{T}=(\mathfrak{K}_{1}^{\epsilon}(\mathfrak{Z}(\cdot)),\mathfrak{K}_{2}^{\epsilon}(\mathfrak{Z}(\cdot))).$$
Assume that $\mathfrak{Z},\tilde{\mathfrak{Z}}\in C_{-\frac{\mu}{\epsilon}}^{+},$ then from (19) we obtain the following upper-bound\begin{align*}||e^{A_{\alpha}t/\epsilon}(X(0)-\tilde{X}(0))||_{1}&\leq e^{-\lambda_{1}t/\epsilon}LipH^{\epsilon}||Y(0)-\tilde{Y}(0)||_{2}\\
&\leq e^{-\lambda_{1}t/\epsilon}LipH^{\epsilon}||\int_{+\infty}^{0}e^{F(-s)}(G(\mathfrak{Z}(s),\theta_{s}^{1}\omega_{1},\theta_{s}^{2}\omega_{2})\\&\indent-G(\tilde{\mathfrak{Z}}(s),\theta_{s}^{1}\omega_{1},\theta_{s}^{2}\omega_{2}))ds||_{2}\\
&\leq e^{-\lambda_{1}t/\epsilon}LipH^{\epsilon}K\int_{0}^{+\infty}e^{-\gamma_{2} s}||\mathfrak{Z}(s)-\tilde{\mathfrak{Z}}(s)||ds,\end{align*}and so,\begin{align*}||\mathfrak{K}_{1}^{\epsilon}(\mathfrak{Z}-\tilde{\mathfrak{Z}})||_{C_{\beta}^{1,+}}&\leq LipH^{\epsilon}\times K||\mathfrak{Z}-\tilde{\mathfrak{Z}}||_{C_{\beta}^{+}}\mathop {\mbox{sup} }\limits_{t\geq 0}\{e^{-(\beta+\frac{\lambda_{1}}{\epsilon})t}\int_{0}^{+\infty}e^{(-\gamma_{2}+\beta)s}ds\}\\&\indent+\frac{K}{\epsilon}||\mathfrak{Z}-\tilde{\mathfrak{Z}}||_{C_{\beta}^{+}}
\mathop {\mbox{sup} }\limits_{t\geq0}\{e^{-\beta t}\int_{0}^{t}e^{-\lambda_{1}(t-s)/\epsilon}ds\}.\end{align*}
This implies that \begin{align}\label{sde}||\mathfrak{K}_{1}^{\epsilon}(\mathfrak{Z}-\tilde{\mathfrak{Z}})||_{C_{\beta}^{1,+}}\leq(\frac{Lip H^{\epsilon}\times K}{\gamma_{2}-\beta}+\frac{K}{\epsilon\beta+\lambda_{1}})||\mathfrak{Z}-\tilde{\mathfrak{Z}}||_{C_{\beta}^{+}}.\end{align}
Similarly, for the operator $\mathfrak{K}_{2}^{\epsilon}$ \begin{align*}||\mathfrak{K}_{2}^{\epsilon}(\mathfrak{Z}-\tilde{\mathfrak{Z}})||_{C_{\beta}^{2,+}}\leq K||\mathfrak{Z}-\tilde{\mathfrak{Z}}||_{C_{\beta}^{+}}\mathop {\mbox{sup} }\limits_{t\geq0}\{e^{-\beta t}\int_{t}^{+\infty}e^{(\gamma_{2})(t-s)}ds\}.\end{align*} Which implies that \begin{align}\label{sde} ||\mathfrak{K}_{2}^{\epsilon}(\mathfrak{Z}-\tilde{\mathfrak{Z}})||_{C_{\beta}^{2,+}}\leq\frac{ K}{\gamma_{2}-\beta}||\mathfrak{Z}-\tilde{\mathfrak{Z}}||_{C_{\beta}^{+}}.\end{align}
From (Theorem 1), we know that
$$LipH^{\epsilon}(\omega,\cdot)\leq\frac{K}{(-\mu+\lambda_{1})[1-K(\frac{1}{-\mu+\lambda_{1}}+\frac{\epsilon}{\epsilon\gamma_{2}+\mu})]}.$$
Now, we get the upper-bound by combining (20) and (21)
$$||\mathfrak{K}^{\epsilon}(\mathfrak{Z}-\tilde{\mathfrak{Z}})||_{C_{-\frac{\mu}{\epsilon}}^{+}}\leq\varrho(K,\lambda_{1},\gamma_{2},\mu,\epsilon))||\mathfrak{Z}-\tilde{\mathfrak{Z}}||_{C_{-\frac{\mu}{\epsilon}}^{+}},$$
where,
$$\varrho(K,\lambda_{1},\gamma_{2},\mu,\epsilon))=\frac{K}{\epsilon\beta+\lambda_{1}}+\frac{\epsilon K}
{\gamma_{2}-\beta}+\frac{ K^{2}}{(-\mu+\lambda_{1})(\gamma_{2}-\beta)[1-K(\frac{1}{-\mu+\lambda_{1}}+\frac{\epsilon}{\epsilon\gamma_{2}+\mu})]}.$$
Take $\beta=-\frac{\mu}{\epsilon}$, then
$$\varrho(K,\lambda_{1},\gamma_{2},\mu,\epsilon))\rightarrow\frac{K}{-\mu+\lambda_{1}}\mbox{ as }\epsilon\rightarrow0.$$
So by (10) there is a sufficiently small constant $\tilde{\epsilon}_{0}>0$ such that
$$\varrho(K,\lambda_{1},\gamma_{2},\mu,\epsilon)<1,\indent \forall\indent 0<\epsilon<\tilde{\epsilon}_{0}.$$
This implies that the operator $\mathfrak{K}^{\epsilon}$ is contractive and also has a unique fixed point $\mathfrak{Z}\in C_{-\frac{\mu}{\epsilon}}^{+}$, this unique fixed point is a unique solution for (18) and satisfies $$(\tilde{U}_{0},\tilde{V}_{0})^{T}=(X(0),Y(0))^{T}+(U_{0},V_{0})^{T}\in M^{\epsilon}(\omega).$$Furthermore,
$$||\mathfrak{Z}||_{C_{-\frac{\mu}{\epsilon}}^{+}}\leq\frac{1}{1-K(\frac{1}{-\mu+\lambda_{1}}+\frac{\epsilon}{\epsilon\gamma_{2}+\mu})}||\mathfrak{Z}(0)||,$$
which means that
$$||\Phi^{\epsilon}(t,\omega,Z_{0})-\Phi^{\epsilon}(t,\omega,\tilde{Z}_{0})||\leq \frac{e^{-\frac{\mu}{\epsilon}t}}{1-K(\frac{1}{-\mu+\lambda_{1}}+\frac{\epsilon}{\epsilon\gamma_{2}+\mu})}||Z_{0}-\tilde{Z}_{0}||,t\geq0.$$
Since the exponential tracking property of $M^{\epsilon}(\omega)$ is obtained, therefore, it completes the proof.  $\Box$\\
With the help of Remark 4 and Remark 5, we get a reduced system based on slow manifold for the original fast-slow stochastic evolutionary system (1)-(2) stated as Theorem 3.\\
\textbf{Theorem 3.} (Slow reduction) Suppose that (H1-H3) hold and $\epsilon>0$ is small enough. For every solution $z^{\epsilon}(t)=(u^{\epsilon}(t),v^{\epsilon}(t))$ to (1)-(2), there is an orbit $\bar{z}^{\epsilon}(t)=(H^{\epsilon}(\omega,\bar{v}^{\epsilon}(t))+\eta^{\epsilon}(\omega_{1}),\bar{v}^{\epsilon}(t)+,\xi(\omega_{2})),$ lying on $\check{M}^{\epsilon}(\omega)$, which satisfies the equation
$$\dot{\bar{v}}^{\epsilon}=J\bar{v}^{\epsilon}+g(H^{\epsilon}(\theta_{t}\omega,\bar{v}^{\epsilon})+\eta^{\epsilon}(\theta_{t}^{1}\omega_{1}),\bar{v}^{\epsilon}+\xi(\theta_{t}^{2}\omega_{2}))+\sigma_{2}W_{t}^{2}$$
such that for every $\omega$ and $t\geq0$,
$$||z^{\epsilon}(t,\omega)-\bar{z}^{\epsilon}(t,\omega)||\leq\frac{e^{-\frac{\mu}{\epsilon}t}}{1-K(\frac{1}{-\mu+\lambda_{1}}+\frac{\epsilon}{\epsilon\gamma_{2}+\mu})}
||z_{0}-\bar{z}_{0}||, \indent t\geq0,$$where$$z_{0}=(u^{\epsilon}(0),v^{\epsilon}(0)),\indent \bar{z}_{0}=(\bar{u}^{\epsilon}(0),\bar{v}^{\epsilon}(0)).$$\end{linenomath*}
\section{Parameter estimation on random slow manifold}
\begin{linenomath*}Assume that original slow equation $(2)$ has an unknown parameter $d \in \mathbb{R}$ and we know its range is denoted by $\mathfrak{\wedge}$, which is a close interval. That is, nonlocal fast-slow system (1)-(2) becomes\begin{align}
&\dot{u}^{\epsilon}=-\frac{1}{\epsilon}(-\Delta)^{\frac{\alpha}{2}}u^{\epsilon}+\frac{1}{\epsilon}f(u^{\epsilon},v^{\epsilon})+\frac{\sigma_{1}}{\sqrt{\epsilon}}\dot{W}_{t}^{1},
\mbox{ in }H_{1},\\&\dot{v}^{\epsilon}=Jv^{\epsilon}+g(u^{\epsilon},v^{\epsilon},d)+\sigma_{2}\dot{W}_{t}^{2},\mbox{ in } H_{2}.\end{align}
In this part, we establish a parameter estimation method for nonlocal stochastic system by using only the observations on slow component, additionaly with error estimation.\end{linenomath*}
\subsection{Parameter estimation based on reduced slow system}\indent

\begin{linenomath*}
\noindent \textbf{Theorem 4.} (Parameter estimation) Suppose that the assumptions (H1-H3) are satisfied and $\epsilon>0$ is small enough. Then nonlocal fast-slow system (22)-(23) has a parameter estimation property in the sense:\\ For every $d$ in nonlocal fast-slow system (22)-(23), there is a $d_{E}^{s}$ in reduced slow system
\begin{align*}\dot{v}^{\epsilon}=Jv^{\epsilon}+g(H^{\epsilon}(\theta_{t}\omega,b^{\epsilon})+\eta^{\epsilon}(\theta_{t}^{1}\omega_{1}),b^{\epsilon}+\xi(\theta_{t}^{2}\omega_{2}),d_{E}^{s})+\sigma_{2}W_{t}^{2},
\mbox{ in }H_{2},\end{align*}
such that
\begin{align*}|d-d_{E}^{s}|&<\frac{1}{G(d,d_{E}^{s})}
\[CL_{g}\mathbb{E}||u_{0}^{\epsilon}-\eta^{\epsilon}(\theta_{t}^{1}\omega_{1})-H^{\epsilon}(\omega,v_{0}^{\epsilon})||_{H_{1}}
\times\frac{\epsilon}{\lambda_{1}-CL_{f}}\\&\indent+\frac{CL_{f}L_{g}}{\lambda_{1}-CL_{f}}
(T\mathbb{E}\int_{0}^{T}||v_{ob}^{\epsilon}(t)-v_{s}^{\epsilon}(t)||_{H_{2}}^{2}dt)^{\frac{1}{2}}+L_{g}(T\mathbb{E}\int_{0}^{T}
||v_{ob}^{\epsilon}(t)\\&\indent-v_{s}^{\epsilon}(t)||_{H_{2}}^{2}dt)^{\frac{1}{2}}+\frac{1}{T}(\mathbb{E}\int_{0}^{T}||v_{ob}^{\epsilon}(t)-v_{s}^{\epsilon}(t)||_{H_{2}}^{2}dt)^{\frac{1}{2}}\].\end{align*}
where $(u^{\epsilon}(0),v^{\epsilon}(0))=(u_{0},v_{0})$, $\mbox{ for }T>0, t\in(0,T),d'\in (d_{E}^{s},d)$ or $d'\in (d,d_{E}^{s}),0<G(d,d_{E}^{s})<\infty,$\\and \\$G(d,d_{E}^{s}):=\mathbb{E}||\int_{0}^{t^{*}}e^{-Jt}\nabla_{d}g(H^{\epsilon}(\theta_{t}\omega,b_{s}^{\epsilon}(t))+\eta^{\epsilon}(\theta_{t}^{1}\omega_{1}),b_{s}^{\epsilon}(t)+\xi(\theta_{t}^{2}\omega_{2}),d')dt||_{H_{2}}.$$
$\\$L_{f}$ and $L_{g}$ are Lipschitz constants of Lipschitz continuous functions $f$ and $g$.\\While $u_{ob}^{\epsilon}(t), v_{ob}^{\epsilon}(t)\mbox{ and } v_{s}^{\epsilon}(t)$ are observations of fast system, slow system and reduced slow system respectively.
\\
 \textbf{ Importance:}
Since slow system has lower dimension than the original nonlocal stochastic system. In computational cost our method offers a benefit of using observations only on slow components. As compared to fast components, it is often much more reliable to observe slow components \cite{ren2015parameter}.\end{linenomath*}

\begin{linenomath*}\noindent\textbf{Proof.} Let observation of $v^{\epsilon}(t)$ is denoted by $v_{ob}^{\epsilon}(t),t\in[0,T],$ corresponding to the original parameter value $d$. We denote the observation of $u^{\epsilon}(t)$ by $u_{ob}^{\epsilon}(t)$ (although, it will not arise in the error estimation described below).\\ We estimate the system parameter $d$ in (22)-(23) using the reduced slow system
\begin{align}\dot{v}^{\epsilon}=Jv^{\epsilon}+g(H^{\epsilon}(\theta_{t}\omega,v^{\epsilon})+\eta^{\epsilon}(\theta_{t}^{1}\omega_{1}),b^{\epsilon}+\xi(\theta_{t}^{2}\omega_{2}),d_{E}^{s})+\sigma_{2}W_{t}^{2},
\mbox{ in }H_{2},\end{align}
for sufficiently small $\epsilon$. Assume that $g(u^{\epsilon},v^{\epsilon},d)$ is Lipschitz continuous with regard to $u^{\epsilon},v^{\epsilon},d$ having Lipschitz constant $L_{g}$. For $v_{s}^{\epsilon}(t)$, which satisfies the reduced system (24) with parameter $d$, represented as $d_{E}^{s}$ just to differentiate from $d$ in (23), and initial data $v_{0}^{\epsilon}$, define an objective function$$
F^{s}(d_{E}^{s})=\mathbb{E}\int_{0}^{T}||v_{ob}^{\epsilon}(t)-v_{s}^{\epsilon}(t)||_{H_{2}}^{2}dt.$$
Take corresponding minimizer $d_{E}^{s}$ as the estimation of original parameter $d$. Now, for this parameter estimation method we provide an error estimation. By using mean value theorem, there is a $t^{*}\in (0,T)$, such that
$$\int_{0}^{T}||v_{ob}^{\epsilon}(t)-v_{s}^{\epsilon}(t)||_{H_{2}}dt=T||v_{ob}^{\epsilon}(t^{*})-v_{s}^{\epsilon}(t^{*})||_{H_{2}},$$
and now by  Cauchy-Schwarz inequality
\begin{equation}||v_{ob}^{\epsilon}(t^{*})-v_{s}^{\epsilon}(t^{*})||_{H_{2}}=\frac{1}{T}\int_{0}^{T}||v_{ob}^{\epsilon}(t)-v_{s}^{\epsilon}(t)||_{H_{2}}dt\leq\frac{1}{T}(T\int_{0}^{T}||v_{ob}^{\epsilon}(t)-v_{s}^{\epsilon}(t)||_{H_{2}}^{2}dt)^{\frac{1}{2}}.\end{equation}
Evaluate the difference between $v_{ob}^{\epsilon}(t)$ and $v_{s}^{\epsilon}(t)$, which satisfies
\begin{align*}\dot{v}_{ob}^{\epsilon}(t)-\dot{v}_{s}^{\epsilon}(t)&=J(v_{ob}^{\epsilon}-v_{s}^{\epsilon})+[g(u_{ob}^{\epsilon},v_{ob}^{\epsilon},d)-g(H^{\epsilon}(\theta_{t}\omega,v_{s}^{\epsilon})\\&\indent+\eta^{\epsilon}(\theta_{t}^{1}\omega_{1}),v_{s}^{\epsilon}+\xi(\theta_{t}^{2}\omega_{2}),d_{E}^{s})],\end{align*} where$$v_{ob}^{\epsilon}(0)-v_{s}^{\epsilon}(0)=0.$$
By using variation of constants formula, it is obtained that\begin{align*}e^{-Jt^{*}}[v_{ob}^{\epsilon}(t^{*})-v_{s}^{\epsilon}(t^{*})]&= \int_{0}^{t*}e^{-Jt}[g(u_{ob}^{\epsilon}(t),v_{ob}^{\epsilon}(t),d)-g(H^{\epsilon}(\theta_{t}\omega,v_{s}^{\epsilon}(t))\\&\indent+\eta^{\epsilon}(\theta_{t}^{1}\omega_{1}),v_{s}^{\epsilon}(t)+\xi(\theta_{t}^{2}\omega_{2}),d)]dt
+\int_{0}^{t^{*}}e^{-Jt}\\&\indent[g(H^{\epsilon}(\theta_{t}\omega,v_{s}^{\epsilon}(t))+\eta^{\epsilon}(\theta_{t}^{1}\omega_{1}),v_{s}^{\epsilon}(t)+\xi(\theta_{t}^{2}\omega_{2}),d)\\&\indent-g(H^{\epsilon}(\theta_{t}\omega,v_{s}^{\epsilon}(t))+\eta^{\epsilon}(\theta_{t}^{1}\omega_{1}),v_{s}^{\epsilon}(t)+\xi(\theta_{t}^{2}\omega_{2}),d_{E}^{s})]dt,
\end{align*}
and then
\begin{align*}&\int_{0}^{t^{*}}e^{-Jt}[g(H^{\epsilon}(\theta_{t}\omega,v_{s}^{\epsilon}(t))+\eta^{\epsilon}(\theta_{t}^{1}\omega_{1}),v_{s}^{\epsilon}(t)+\xi(\theta_{t}^{2}\omega_{2}),d)-g(H^{\epsilon}(\theta_{t}\omega,v_{s}^{\epsilon}(t))\end{align*}\begin{align*}&+\eta^{\epsilon}(\theta_{t}^{1}\omega_{1}),v_{s}^{\epsilon}(t)+\xi(\theta_{t}^{2}\omega_{2}),d_{E}^{s})]dt=-
\int_{0}^{t*}e^{-Jt}[g(u_{ob}^{\epsilon}(t),v_{ob}^{\epsilon}(t),d)\end{align*}\begin{align}-g(H^{\epsilon}(\theta_{t}\omega,v_{s}^{\epsilon}(t))+\eta^{\epsilon}(\theta_{t}^{1}\omega_{1}),v_{s}^{\epsilon}(t)+\xi(\theta_{t}^{2}\omega_{2}),d)]dt+
e^{-Jt^{*}}[v_{ob}^{\epsilon}(t^{*})-v_{s}^{\epsilon}(t^{*})].\end{align}
Taking norm on both sides, and by using the mean value theorem, there is a $k\in (0,1),$ and $d'=d+k(d_{E}^{s}+d)$ on left hand side of (26). Together by using hypothesis $H1$, Lipschitz continuity of $g$, jointly with Cauchy-Schwartz inequality on right hand side of (26) lead to the upper-bound
\begin{align*}&|d-d_{E}^{s}|\cdot||\int_{0}^{t^{*}}e^{-Jt}\nabla_{d}g(H^{\epsilon}(\theta_{t}\omega,v_{s}^{\epsilon}(t))+\eta^{\epsilon}(\theta_{t}^{1}\omega_{1}),b_{s}^{\epsilon}(t)+\xi(\theta_{t}^{2}\omega_{2}),d')dt||_{H_{2}}
\end{align*}\begin{align*}\indent&\leq\int_{0}^{t^{*}}e^{-\gamma_{2}t}L_{g}[||u_{ob}^{\epsilon}(t)-H^{\epsilon}(\theta_{t}\omega,v_{s}^{\epsilon}(t))-\eta^{\epsilon}(\theta_{t}^{1}\omega_{1}))||_{H_{1}}
\end{align*}\begin{align*}&\noindent+||v_{ob}^{\epsilon}(t)-v_{s}^{\epsilon}(t)||_{H_{2}}]dt+e^{-\gamma_{2}t}||v_{ob}^{\epsilon}(t^{*})-v_{s}^{\epsilon}(t^{*})||_{H_{2}}
\end{align*}\begin{align*}<L_{g}\int_{0}^{T}||u_{ob}^{\epsilon}(t)-H^{\epsilon}(\theta_{t}\omega,v_{s}^{\epsilon}(t))-\eta^{\epsilon}(\theta_{t}^{1}\omega_{1}))||_{H_{1}}\end{align*}\begin{align}\indent+L_{g}(T\int_{0}^{T}
||v_{ob}^{\epsilon}(t)-v_{s}^{\epsilon}(t)||_{H_{2}}^{2}dt)^{\frac{1}{2}}+||v_{ob}^{\epsilon}(t^{*})-v_{s}^{\epsilon}(t^{*})||_{H_{2}}.\end{align}
As stated, the random transformation converts stochastic evolutionary system to random system. Thus for $(u_{ob}^{\epsilon}(t),v_{ob}^{\epsilon}(t))$ satisfying the stochastic system (1)-(2), the corresponding transformation $(u_{ob}^{\epsilon}(t)-\eta^{\epsilon}(\theta_{t}^{1}\omega_{1}),v_{ob}^{\epsilon}(t)-\xi(\theta_{t}^{2}\omega_{2}))$ satisfies the random system (8)-(9) with an unknown parameter value $d$ and initial value
$(u_{0}^{\epsilon}-\eta^{\epsilon}(\theta_{t}^{1}\omega_{1}),v_{0}^{\epsilon}-\xi(\theta_{t}^{2}\omega_{2}))$. As we are familiar that special paths of corresponding system are consist in the slow manifold. Thus $(H^{\epsilon}(\theta_{t}\omega,v_{s}^{\epsilon}(t)),v_{s}^{\epsilon}(t))$ also satisfies (8)-(9) with parameter denoted as $d_{E}^{s}$ and initial value $(H^{\epsilon}(\theta_{t}\omega,v_{0}^{\epsilon}),v_{0}^{\epsilon})$.\\
From random dynamical system (8), the difference between $H^{\epsilon}(\theta_{t}\omega,v_{s}^{\epsilon}(t))+\eta^{\epsilon}(\theta_{t}^{1}\omega_{1})$ and $u_{ob}^{\epsilon}(t)$ satisfies,\\ \\
$\frac{d}{dt}[u_{ob}^{\epsilon}(t)-\eta^{\epsilon}(\theta_{t}^{1}\omega_{1})-H^{\epsilon}(\theta_{t}\omega,v_{s}^{\epsilon}(t))]
=\frac{1}{\epsilon}A_{\alpha}[u_{ob}^{\epsilon}(t)-\eta^{\epsilon}(\theta_{t}^{1}\omega_{1})$\\ \\ $-H^{\epsilon}(\theta_{t}\omega,v_{s}^{\epsilon}(t))]
+\frac{1}{\epsilon}[f(u_{ob}^{\epsilon}(t),v_{ob}^{\epsilon}(t))-f(H^{\epsilon}(\theta_{t}\omega,v_{s}^{\epsilon}(t))+\eta^{\epsilon}(\theta_{t}^{1}\omega_{1}),v_{s}^{\epsilon}(t)+\xi(\theta_{t}^{2}\omega_{2}))].$
\\ \\By using variation of constants formula, we obtain
\begin{align*}&||u_{ob}^{\epsilon}(t)-\eta^{\epsilon}(\theta_{t}^{1}\omega_{1})-H^{\epsilon}(\theta_{t}\omega,v_{s}^{\epsilon}(t))||_{H_{1}}
\\&\indent=||e^{\frac{A_{\alpha}}{\epsilon}t}[u_{0}^{\epsilon}-\eta^{\epsilon}(\theta_{t}^{1}\omega_{1})-H^{\epsilon}(\theta_{t}\omega,v_{0}^{\epsilon})]+\frac{1}{\epsilon}\int_{0}^{t}
e^{\frac{A_{\alpha}}{\epsilon}(t-s)}[f(u_{ob}^{\epsilon}(s),v_{ob}^{\epsilon}(s))
\\&\indent \indent-f(H^{\epsilon}(\theta_{t}\omega,v_{s}^{\epsilon}(s))+\eta^{\epsilon}(\theta_{t}^{1}\omega_{1}),v_{s}^{\epsilon}(s)+\xi(\theta_{t}^{2}\omega_{2}))]ds||_{H_{1}}
\\&\indent\leq Ce^{-\frac{\lambda_{1}}{\epsilon}t}||u_{0}^{\epsilon}-\eta^{\epsilon}(\theta_{t}^{1}\omega_{1})-H^{\epsilon}(\theta_{t}\omega,v_{0}^{\epsilon})||_{H_{1}}
+\frac{CL_{f}}{\epsilon}\int_{0}^{t}e^{\frac{\lambda_{1}}{\epsilon}(s-t)}[||u_{ob}^{\epsilon}(s)
\\&\indent\indent-(H^{\epsilon}(\theta_{t}\omega,v_{s}^{\epsilon}(s))-\eta^{\epsilon}(\theta_{t}^{1}\omega_{1}))||_{H_{1}}+||v_{ob}^{\epsilon}(s))-v_{s}^{\epsilon}(s)||_{H_{2}}]ds.\end{align*}
This implies that\begin{align*}&e^{\frac{\lambda_{1}}{\epsilon}t}||u_{ob}^{\epsilon}(t)-\eta^{\epsilon}(\theta_{t}^{1}\omega_{1})-H^{\epsilon}(\theta_{t}\omega,v_{s}^{\epsilon}(t))||_{H_{1}}
\\&\indent\leq C||u_{0}^{\epsilon}-\eta^{\epsilon}(\theta_{t}^{1}\omega_{1})-H^{\epsilon}(\theta_{t}\omega,v_{0}^{\epsilon})||_{H_{1}}
+\frac{CL_{f}}{\epsilon}\int_{0}^{t}e^{\frac{\lambda_{1}}{\epsilon}s}||u_{ob}^{\epsilon}(s)
\\&\indent \indent-(H^{\epsilon}(\theta_{t}\omega,v_{s}^{\epsilon}(s))-\eta^{\epsilon}(\theta_{t}^{1}\omega_{1}))||_{H_{1}}ds+\frac{CL_{f}}{\epsilon}\int_{0}^{t}e^{\frac{\lambda_{1}}{\epsilon}s}||v_{ob}^{\epsilon}(s))
\\&\indent \indent-v_{s}^{\epsilon}(s)||_{H_{2}}ds.\end{align*}
By Gronwall's inequality, we obtained
\begin{align*}&e^{\frac{\lambda_{1}}{\epsilon}t}||u_{ob}^{\epsilon}(t)-\eta^{\epsilon}(\theta_{t}^{1}\omega_{1})-H^{\epsilon}(\theta_{t}\omega,v_{s}^{\epsilon}(t))||_{H_{1}}
\\&\indent\leq C||u_{0}^{\epsilon}-\eta^{\epsilon}(\theta_{t}^{1}\omega_{1})-H^{\epsilon}(\omega,v_{0}^{\epsilon})||_{H_{1}}e^{\frac{CL_{f}}{\epsilon}t}
+\frac{CL_{f}}{\epsilon}\int_{0}^{t}e^{\frac{\lambda_{1}}{\epsilon}s}||v_{ob}^{\epsilon}(s)
\\&\indent\indent-v_{s}^{\epsilon}(s)||_{H_{2}}e^{\frac{CL_{f}}{\epsilon}(t-s)}ds,\end{align*}
and thus\begin{align*}&||u_{ob}^{\epsilon}(t)-\eta^{\epsilon}(\theta_{t}^{1}\omega_{1})-H^{\epsilon}(\theta_{t}\omega,v_{s}^{\epsilon}(t))||_{H_{1}}
\\&\indent\indent\leq C||u_{0}^{\epsilon}-\eta^{\epsilon}(\theta_{t}^{1}\omega_{1})-H^{\epsilon}(\omega,v_{0}^{\epsilon})||_{H_{1}}e^{-\frac{\lambda_{1}-CL_{f}}{\epsilon}t}
+\frac{CL_{f}}{\epsilon}\int_{0}^{t}||v_{ob}^{\epsilon}(s)\\&\indent\indent \indent
-v_{s}^{\epsilon}(s)||_{H_{2}}e^{-\frac{\lambda_{1}-CL_{f}}{\epsilon}(t-s)}ds.\end{align*}
Now, we conclude by exchanging the order of integrals,
\begin{align*}&\int_{0}^{T}\int_{0}^{t}||v_{ob}^{\epsilon}(s)-v_{s}^{\epsilon}(s)||_{H_{2}}e^{-\frac{\lambda_{1}-CL_{f}}{\epsilon}(t-s)}dsdt\\&\indent=\int_{0}^{T}\int_{s}^{T}||v_{ob}^{\epsilon}(s) -v_{s}^{\epsilon}(s)||_{H_{2}}e^{-\frac{\lambda_{1}-CL_{f}}{\epsilon}(t-s)}dtds \\&\indent<\frac{\epsilon}{\lambda_{1}-C L_{f}}\int_{0}^{T}||v_{ob}^{\epsilon}(s)-v_{s}^{\epsilon}(s)||_{H_{2}} ds.\end{align*}
Thus by Cauchy-Schwartz inequality
\begin{align*}
&\int_{0}^{T}||u_{ob}^{\epsilon}(t)-\eta^{\epsilon}(\theta_{t}^{1}\omega_{1})-H^{\epsilon}(\theta_{t}\omega,v_{s}^{\epsilon}(t))||_{H_{1}}dt
<C||u_{0}^{\epsilon}-\eta^{\epsilon}(\theta_{t}^{1}\omega_{1})\end{align*}\begin{align*}&-H^{\epsilon}(\omega,v_{0}^{\epsilon})||_{H_{1}}
\times\frac{\epsilon(1-e^{-\frac{\lambda_{1}-CL_{f}}{\epsilon}T})}{\lambda_{1}-CL_{f}}
+\frac{CL_{f}}{\lambda_{1}-CL_{f}}\int_{0}^{T}||v_{ob}^{\epsilon}(s)\end{align*}\begin{align*}&-v_{s}^{\epsilon}(s)||_{H_{2}}ds
<C||u_{0}^{\epsilon}-\eta^{\epsilon}(\theta_{t}^{1}\omega_{1})-H^{\epsilon}(\omega,v_{0}^{\epsilon})||_{H_{1}}\frac{\epsilon}{\lambda_{1}-CL_{f}}\end{align*}\begin{align}
+\frac{CL_{f}}{\lambda_{1}-CL_{f}}\left(T\int_{0}^{T}||v_{ob}^{\epsilon}(t)-v_{s}^{\epsilon}(t)||_{H_{2}}^{2}dt\right)^{\frac{1}{2}}.\end{align}
Inserting (25) and (28) into (27), and taking expectation on both sides, also by using Cauchy-Schwartz inequality $\mathbb{E}|XY|\leq\mathbb{E}(|X|^{2})^{\frac{1}{2}}\mathbb{E}(|Y|^{2})^{\frac{1}{2}}$, we obtain the following estimation
\begin{align*}&|d-d_{E}^{s}|\cdot\mathbb{E}||\int_{0}^{t^{*}}e^{-Jt}\nabla_{d}g(H^{\epsilon}(\theta_{t}\omega,v_{s}^{\epsilon}(t))+\eta^{\epsilon}(\theta_{t}^{1}\omega_{1}),b_{s}^{\epsilon}(t)+\xi(\theta_{t}^{2}\omega_{2}),d')dt||_{H_{2}}
\end{align*}\begin{align*}&<CL_{g}\mathbb{E}||u_{0}^{\epsilon}-\eta^{\epsilon}(\theta_{t}^{1}\omega_{1})-H^{\epsilon}(\omega,v_{0}^{\epsilon})||_{H_{1}}
\times\frac{\epsilon}{\lambda_{1}-CL_{f}}+\frac{CL_{f}L_{g}}{\lambda_{1}-CL_{f}}
\end{align*}\begin{align*}&\times\left(T\mathbb{E}\int_{0}^{T}||v_{ob}^{\epsilon}(t)-v_{s}^{\epsilon}(t)||_{H_{2}}^{2}dt\right)^{\frac{1}{2}}+L_{g}\left(T\mathbb{E}\int_{0}^{T}
||v_{ob}^{\epsilon}(t)-v_{s}^{\epsilon}(t)||_{H_{2}}^{2}dt\right)^{\frac{1}{2}}\end{align*}\begin{align}&+\frac{1}{T}\left(\mathbb{E}\int_{0}^{T}||v_{ob}^{\epsilon}(t)-v_{s}^{\epsilon}(t)||_{H_{2}}^{2}dt\right)^{\frac{1}{2}}.\end{align}
Consider integral term, which is on left hand side of (29), and set
\begin{equation*}\label{sde}
G(d,d_{E}^{s}):=\mathbb{E}||\int_{0}^{t^{*}}e^{-Jt}\nabla_{d}g(H^{\epsilon}(\theta_{t}\omega,v_{s}^{\epsilon}(t))+\eta^{\epsilon}(\theta_{t}^{1}\omega_{1}),v_{s}^{\epsilon}(t)+\xi(\theta_{t}^{2}\omega_{2}),d')dt||_{H_{2}},
\end{equation*}
where $T>0, t^{*}\in(0,T),d'\in (d_{E}^{s},d)$ or $d'\in (d,d_{E}^{s}).$\\
Assume the function $G(d,d_{E}^{s})$ satisfies $$0<G(d,d_{E}^{s})<\infty.$$
The estimation indicates that $|d-d_{E}^{s}|$ can be controlled by error due to slow reduction and observation error (objective function), i.e. $\mathcal{O}(\epsilon)$ and $\mathcal{O}((F^{s}(d_{E}^{s}))^{\frac{1}{2}}).$ $\Box$\\
\indent If we can expand $H^{\epsilon}(\theta_{t}\omega,v^{\epsilon}(t))$ for sufficiently small $\epsilon$ up to order $\mathcal{O}(\epsilon^{2})$ like \cite{bai2017slow} as follows
\begin{equation*}\label{sde}
\check{H}^{\epsilon}(\theta_{t}\omega,B_{0})=H^{0}(\omega,V_{0})+\epsilon( H^{1}(\omega,V_{0})+\eta^{\epsilon}(\theta_{t}^{1}\omega_{1}))+\mathcal{O}(\epsilon^{2}).\end{equation*}
Then approximate reduced slow system is
\begin{align}\dot{v}^{\epsilon}=Jv^{\epsilon}+g(\check{H}^{\epsilon}(\theta_{t}\omega,v^{\epsilon})+\eta^{\epsilon}(\theta_{t}^{1}\omega_{1}),b^{\epsilon}+\xi(\theta_{t}^{2}\omega_{2}))+\sigma_{2}W_{t}^{2},
\mbox{ in }H_{2},\end{align}
for the estimation. Since
\begin{align*}
&||u_{ob}^{\epsilon}(t)-\check{H}^{\epsilon}(\theta_{t}\omega,v_{s}^{\epsilon}(t))-\eta^{\epsilon}(\theta_{t}^{1}\omega_{1})||_{H_{1}}
\\&\indent\leq||u_{ob}^{\epsilon}(t)-\eta^{\epsilon}(\theta_{t}^{1}\omega_{1})-H^{\epsilon}(\theta_{t}\omega,v_{s}^{\epsilon}(t))||_{H_{1}}
+||H^{\epsilon}(\theta_{t}\omega,v_{s}^{\epsilon}(t))\\&\indent \indent -\check{H}^{\epsilon}(\theta_{t}\omega,v_{s}^{\epsilon}(t))||_{H_{1}}.\end{align*}
\indent Estimation using (30) is thus bounded by $\mathcal{O}(\epsilon)$ and $\mathcal{O}((F^{s}(d_{E}^{s}))^{\frac{1}{2}}).$\\
\end{linenomath*}

\section{Examples}
\begin{linenomath}To observe the results of existence of slow manifolds, which are obtained in previous section, let us see a few examples.\\
\textbf{Example 1.} Take a system
\begin{align}
&u_{t}=\frac{1}{\epsilon}A_{\alpha}u+\frac{1}{\epsilon}f(u,v)+\frac{1}{\sqrt{\epsilon}}\dot{W_{t}^{1}},\\
&v_{t}=g(u,v)+\dot{W_{t}^{2}},\\ &u|(-1,1)^{c}=0, \mbox{ and } v|(-1,1)^{c}=0.\end{align}
Where $f:L^{2}(-1,1)\times H_{2}\rightarrow L^{2}(-1,1)$ and $g:L^{2}(-1,1)\times H_{2}\rightarrow H_{2}$ are Lipschitz functions with a positive Lipschitz constant $K$:$$|f(u_{1},v_{1})-f(u_{2},v_{2})|\leq K|u_{1}-u_{2}|_{L^{2}(-1,1)}+|v_{1}-v_{2}|_{H_{2}},$$
$$|g(u_{1},v_{1})-g(u_{2},v_{2})|\leq K|u_{1}-u_{2}|_{L^{2}(-1,1)}+|v_{1}-v_{2}|_{H_{2}},$$
for all $(u_{1},v_{1}),(u_{2},v_{2})\in L^{2}(-1,1)\times H_{2}$. Thus hypothesis $(H2)$ holds. As the nonlocal Laplacian operator $A_{\alpha}=-(-\Delta)^{\frac{\alpha}{2}}$ is a generator of $C_{0}$-semigroup $\{e^{A_{\alpha}t}:t\geq0\}$ on $H_{1}=L^{2}(-1,1)$ and also a sectorial operator, which satisfy the upper-bound:
$$||e^{A_{\alpha}t}||_{L^{2}(-1,1)}\leq Ce^{-\lambda_{1}t},$$
with $C>0$, and \\$$||e^{Jt}||_{H_{2}}\leq e^{-\gamma t},\mbox{ for }t\leq0$$with $\gamma=0$. Therefore the system (31)-(32) has a random manifold $M^{\epsilon}(\omega)=\{(H^{\epsilon}(\omega,V),V)^{T}:V\in H_{2}\}$ with exponential tracking property if $K<\frac{\gamma\lambda_{1}}{2\lambda_{1}+\gamma}$ and $\epsilon>0$ is small enough.\\
\indent The above system may be model of a certain biological process \cite{fitzhugh1961impulses,nagumo1962active,cronin1987mathematical}.\\
\end{linenomath}
In Example2, we illustrate our analytical results via numerical simulation.
\noindent \textbf{Example 2.}
\begin{linenomath}Take a fast-slow stochastic system
\begin{align}
&\dot{u}^{\epsilon}=\frac{1}{\epsilon}A_{\alpha}u^{\epsilon}+\frac{0.01}{\epsilon}(\sqrt{v^{\epsilon^{2}}+5}-\sqrt{5})+\frac{\sigma_{1}}{\sqrt{\epsilon}}\dot{W}_{t}^{1}, \mbox{ in } H_{1}=L^{2}(-1,1),\\
&\dot{v}^{\epsilon}=-v^{\epsilon}+(0.01\times a)\sin\int_{-1}^{1}u^{\epsilon}dx+\sigma_{2}\dot{W}_{t}^{2}, \mbox{ in } H_{2}=\mathbb{R},\end{align}
where $u^{\epsilon}$ is fast variable, $v^{\epsilon}$ is slow variable, $a$ is a positive real unknown parameter. While $\dot{W}_{t}^{1}$ and $\dot{W}_{t}^{2}$ are time derivatives of time dependent Weiner processes $W^{1}$ and $W^{2}$ respectively. Nonlinear functions $f=0.01(\sqrt{v^{\epsilon^{2}}+5}-\sqrt{5})$ and $g=(0.01\times a)\sin\int_{-1}^{1}u^{\epsilon}dx$ are Lipschitz continuous with Lipschitz constants $L_{f}=0.01$ and $L_{g}=0.01\times a$ respectively. Stochastic system (38)-(39) is transformed into the random system
\begin{align}
&\dot{U}^{\epsilon}=\frac{1}{\epsilon}A_{\alpha}U^{\epsilon}+\frac{0.01}{\epsilon}(\sqrt{(V^{\epsilon}+\xi(\theta_{t}^{2}\omega_{2}))^{2}+5}-\sqrt{5})\\
&\dot{V}^{\epsilon}=-V^{\epsilon}+(0.01\times a)\sin\left(\int_{-1}^{1}[U^{\epsilon}+\eta^{\epsilon}(\theta_{t}^{1}\omega_{1})]dx\right).\end{align}
Random dynamical system $(40)-(41)$ has a slow invariant manifold, for sufficiently small $\epsilon>0$ $$M^{\epsilon}(\omega)=\{(H^{\epsilon}(\omega,V),V)^{T}:V\in \mathbb{R}),$$
where
$$H^{\epsilon}(\omega,V)=\frac{0.01}{\epsilon}\int_{-\infty}^{0}e^{-A_{\alpha}s/\epsilon}(\sqrt{\left(V^{\epsilon}+\xi(\theta_{t}^{2}\omega_{2})\right)^{2}+5}-\sqrt{5})ds.$$
Based on \cite{bai2017slow,ren2015approximation}, we plot the graph of slow manifolds approximately to order $\mathcal{O}(\epsilon)$. Let $t=\tau\epsilon$\\
$$U^{\epsilon}(t)=U^{\epsilon}(\tau\epsilon)=U_{0}(\tau)+\epsilon U_{1}(\tau)+\epsilon^{2}U_{2}(\tau)+\cdot\cdot\cdot,$$
with initial value
$U^{\epsilon}(0)=H^{\epsilon}(\omega,V_{0})=H^{0}+\epsilon H^{1}+\cdot\cdot\cdot,$ and $V^{\epsilon}(0)=V_{0}.$
$$H^{0}(\omega,V_{0})=\frac{0.01}{\epsilon}\int_{-\infty}^{0}e^{-A_{\alpha}s}(\sqrt{(V_{0})^{2}+5}-\sqrt{5})ds.$$
Then, approximation of $\breve{H}^{\epsilon}(\omega,V_{0})$ upto error $\mathcal{O}(\epsilon)$ is
$$\breve{H}^{\epsilon}(\omega,V_{0})=H^{0}(\omega,V_{0})+\eta^{\epsilon}(\omega).$$
Finally, the approximate reduced slow system is
$$\dot{\bar{v}}^{\epsilon}=-\bar{v}^{\epsilon}+(0.01\times a)\sin\left(\int_{-1}^{1}[H^{0}(\theta_{t}\omega,\bar{v_{0}})+\eta^{\epsilon}(\omega)]\right)dx.$$
We conducted the following numerical simulations to illustrate slow manifold, exponential attracting property and parameter estimation based on slow system.\end{linenomath}
\begin{figure}[http]
  \centering
  \begin{minipage}[b]{0.4\textwidth}
  \includegraphics[width=6cm,height=4cm]{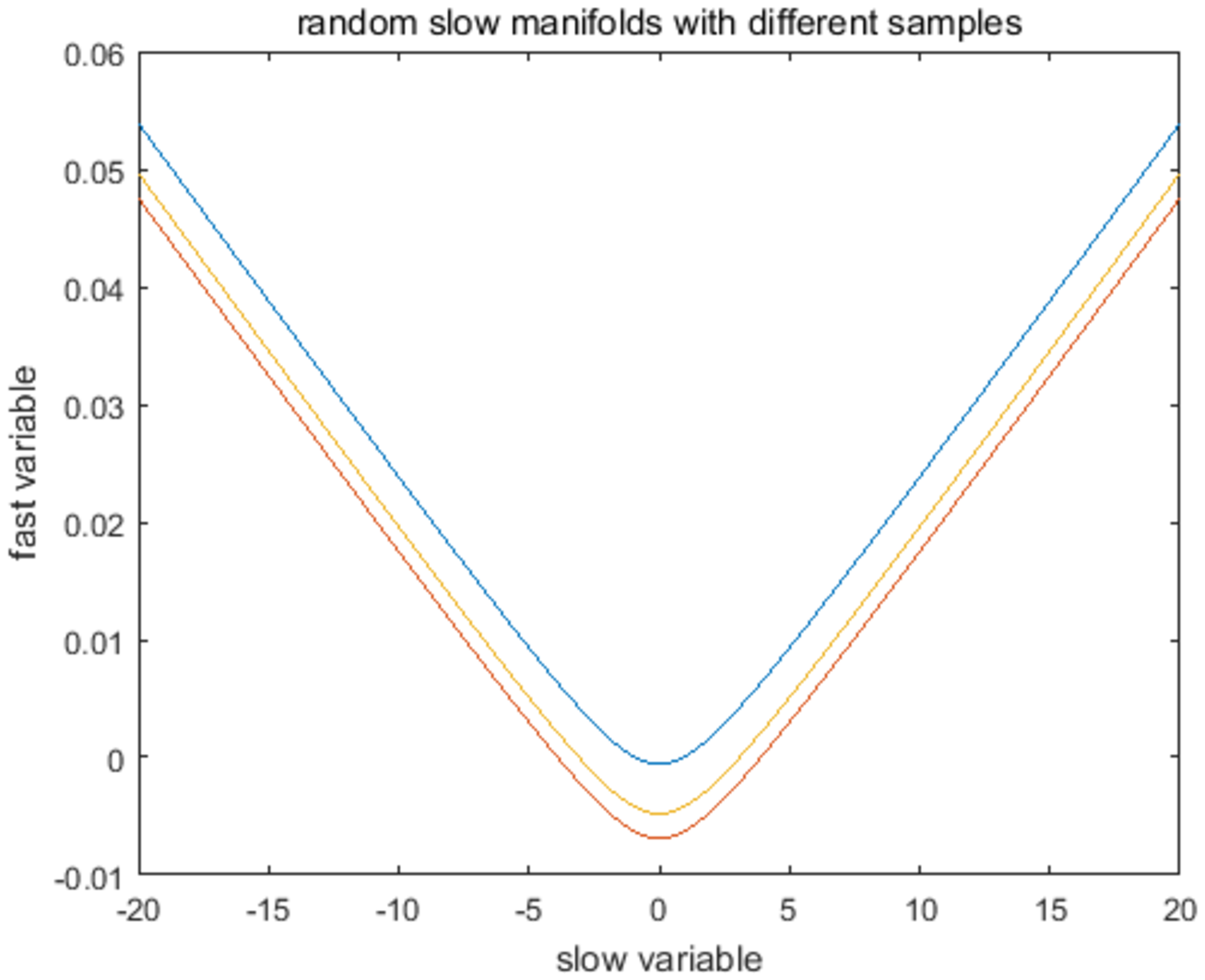}
   \label{fig:graph1}
  \end{minipage}
  \hfill
  \begin{minipage}[b]{0.4\textwidth}
   \includegraphics[width=6cm,height=4cm]{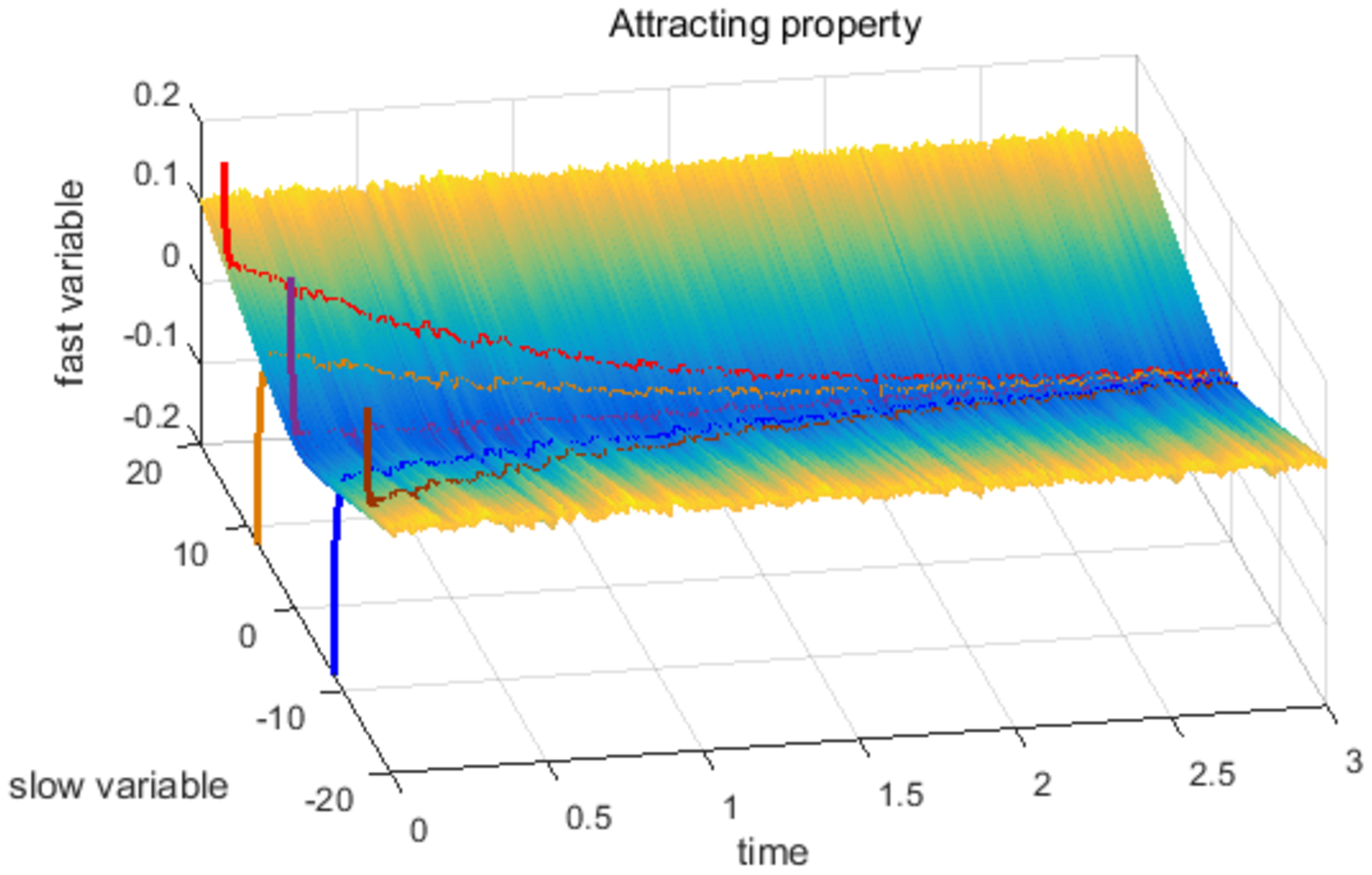}
   \label{fig:graph1}

  \end{minipage}
  \caption{ Random slow manifold for different samples. Each curve stands for a sample
   (left) and exponential tracking property in the system (right): $\sigma_{1}=\sigma_{2}=0.1,\alpha=1.2$ and $\epsilon=0.01$.}
\end{figure}
\newpage
\begin{figure}[http]
  \centering
  \begin{minipage}[b]{0.4\textwidth}
  \includegraphics[width=6cm,height=4cm]{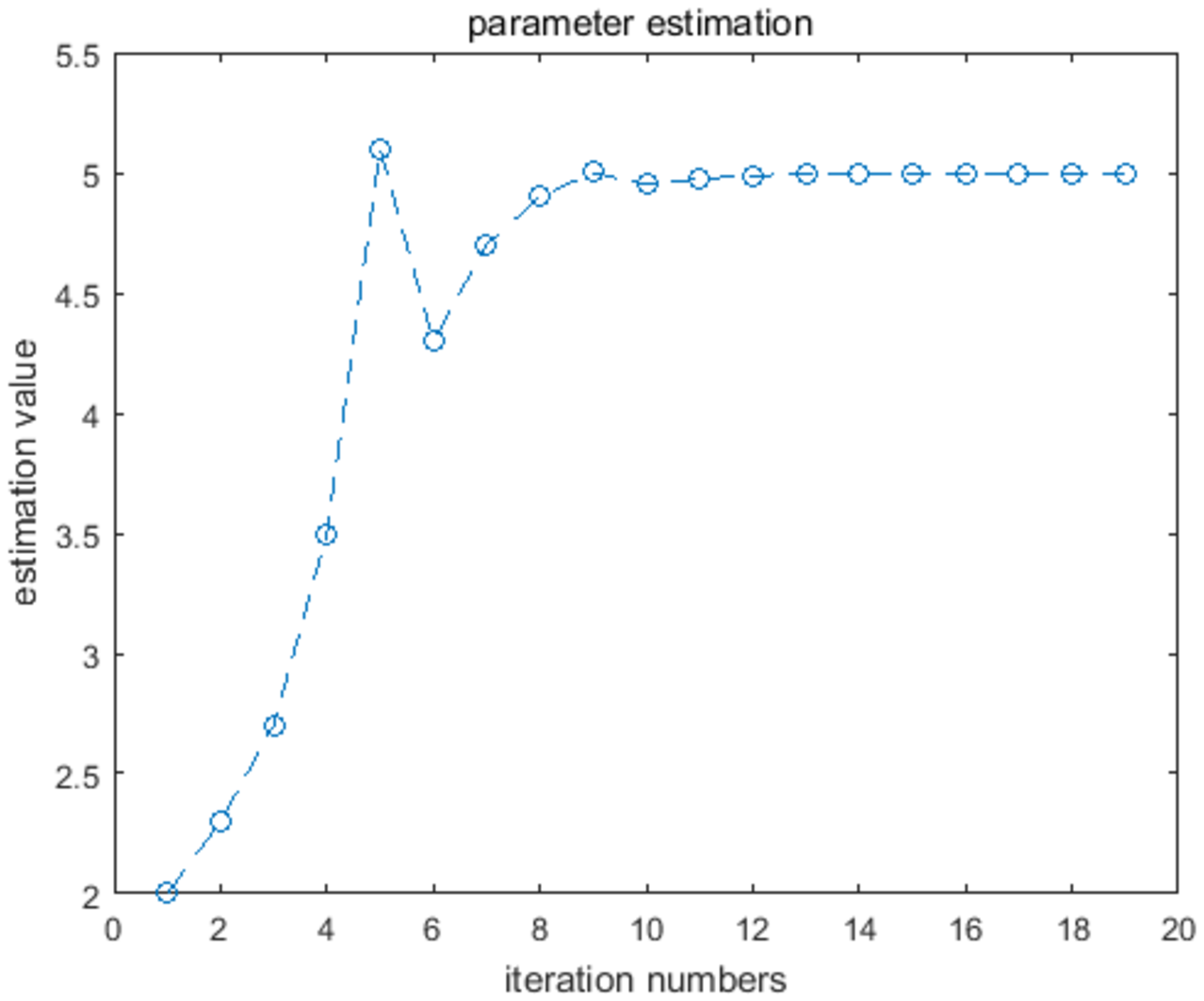}
   \label{fig:graph1}
  \end{minipage}
  \hfill
  \begin{minipage}[b]{0.4\textwidth}
  \includegraphics[width=6cm,height=4cm]{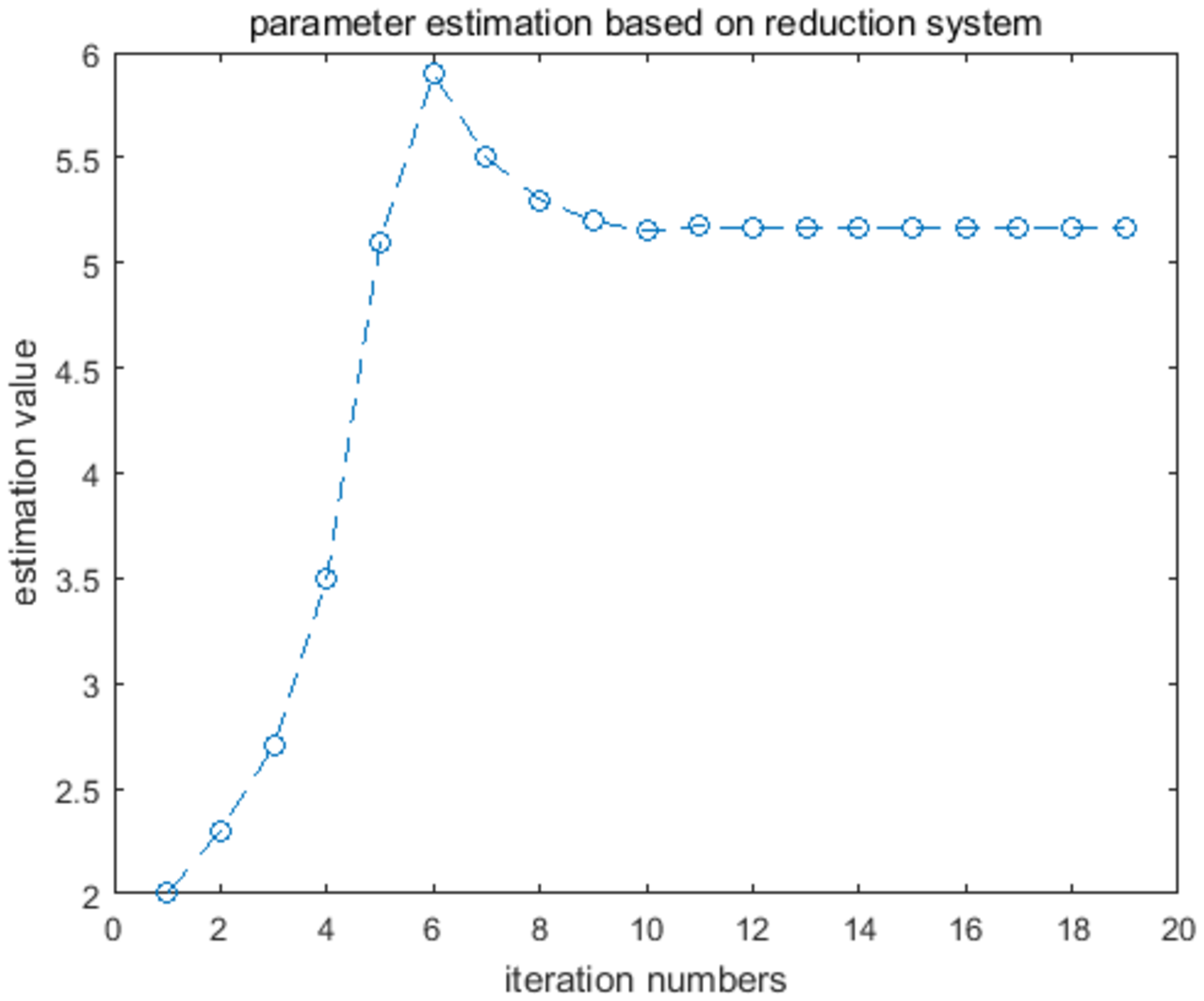}
   \label{fig:graph1}

  \end{minipage}
  \caption{Parameter estimation for the original fast-slow nonlocal stochastic system (left) and parameter estimation for the reduced slow system (right): $\sigma_{1}=\sigma_{2}=0.1,\alpha=1.2$ and $\epsilon=0.01$.}
\end{figure}

\bibliographystyle{line}
\bibliography{JAMS-paper}
\end{document}